
\input amstex
\input amssym.tex
\input amssym.def

\documentstyle{amsppt}
\magnification = \magstep1
\font\bigbf=cmbx10 scaled \magstep2

\vsize=8.55truein

\baselineskip=14pt
\overfullrule=0pt



\define\cCu{{{\Cal Cu}}}


\define\rM{{{\text M}}}


\define\K{{\frak K}}
\define\smallsim{^{\frak{\sim}}}


\define\cC{{{\Cal C}}}

\define\cO{{{\Cal O}}}

\define\IC{{{\Bbb C}}}

\vglue .5in

\centerline{\bigbf The Cuntz semigroup as an invariant for 
C*-algebras}

\vskip .25in

\centerline{Kristofer T. Coward, George A. Elliott, and
Cristian Ivanescu}  

\vskip .3truein 

\midinsert\narrower\narrower

\noindent{\bf Abstract.}
A category is described to which the Cuntz semigroup
belongs and as a functor into which it preserves
inductive limits.
\endinsert

\vskip .3in
{\bf 1.}\
Recently, Toms in [26] used the refinement of the
invariant K$_0$ introduced by Cuntz almost
thirty years ago in [4] to show that certain
C*-algebras are not isomorphic. Anticipating
the possible use of this invariant
to establish isomorphism, we take the liberty
of reporting some observations concerning it. (In particular,
we present what might be viewed as  an embryonic isomorphism
theorem.)

One of the first things that might be noted
in connection with this invariant, which considers,
instead of the finitely generated projective modules
over a given C*-algebra, the larger class of modules consisting of
the countably generated
Hilbert C*-modules (see [12], [15], and [21]; see also [10])  is that, whereas the
equivalence relation between finitely generated projective modules
would appear to be inevitable, namely, just
isomorphism, in the wider setting of Hilbert C*-modules
it is no longer quite so clear what the equivalence
relation should be. While it is tempting just to
choose isomorphism again, one should note that,
even in the stably finite case (which is perhaps
the case that this invariant is of most interest),
whereas the isomorphism classes of algebraically finitely 
generated Hilbert C*-modules (which are of course also algebraically
projective, and up to isomorphism exhaust the
finitely generated projective modules) form an
ordered set with respect to inclusion (in other words, if each
of two such modules is isomorphic to a submodule
of the other, then they must be isomorphic---indeed, any
two such isomorphisms, from each of

\noindent
\underbar{\hskip 2truein}

The research of the second and third authors was supported by 
grants from the Natural Sciences and Engineering
Research Council of Canada.

\vfill\eject

\noindent
two such modules
onto a submodule of the other, must, by stable finiteness,
already be surjective---and so constitute an isomorphism of
the two given modules), this  
would seem to be almost completely
open for isomorphism classes of Hilbert C*-modules.
(Note that we are referring here to the stably finite
case; interestingly, in the purely infinite simple case,
the countably generated
Hilbert C*-modules which are not finitely
generated as modules are all isomorphic---see 4.1.3 of [24]---and
in particular the inclusion relation for isomorphism classes is
antisymmetric! In the case of real rank zero and stable
rank one, which is a very small subset of the stably
finite case, the proof that the Murray-von Neumann equivalence class of a
multiplier projection of an AF algebra is determined by the set of
equivalence classes of projections in the C*-algebra that it majorizes---see [5]---which uses
cancellation, known (see 6.5.1 and 6.5.2 of [1]) to be equivalent to stable rank
one in the
presence of real rank zero---shows that the isomorphism
classes of countably generated
Hilbert C*-modules are determined by
the isomorphism classes of the algebraically finitely
generated ones they include, and it follows again
immediately that the inclusion relation 
on isomorphism classes is antisymmetric.)
(In Theorem 3 below we shall show that this holds without assuming real rank 
zero---provided that stable rank one is still assumed.)
(Note that antisymmetry fails for inclusion of isomorphism
classes of finitely generated projective modules in the
Cuntz algebra $\cO_n$ for $n\ge 3$.) Nevertheless, one would obtain an ordered semigroup
by just antisymmetrizing the inclusion relation on isomorphism classes. 

Cuntz's choice of equivalence relation is a
weaker one again, of an approximative nature (although
there still seems to be no indication whether it is
different from the one just described---i.e., two-way
inclusion of isomorphism classes---or, in the stably
finite case, from just isomorphism for that matter).
It should be emphasized that the present formulation
of Cuntz's invariant, in terms of (countably
generated) Hilbert C*-modules (it could be done
alternatively---but slightly less conveniently---in particular
one would need Brown's theorem to construct suprema!---in
terms of (singly generated) hereditary
sub-C*-algebras of the stabilization of the
C*-algebra) is important for our determination of
the abstract structure of this invariant (see below).
Accordingly, let us describe Cuntz's equivalence relation
in the setting of Hilbert C*-modules.

Cuntz's equivalence relation, like the one considered above, is
based on a pre-order
relation, compatible with direct sum of Hilbert C*-modules
and so again giving rise to an ordered abelian
semigroup (possibly the same one!). For
our purposes, it seems
most appropriate to describe it in terms of the
notion of compact inclusion: let us say that a closed 
submodule of a given Hilbert C*-module over a given 
C*-algebra---let us call this a subobject---is compactly
contained as a subobject---and denote this relation by
$\subset\subset$---if
there is a compact self-adjoint endomorphism of the larger
Hilbert C*-module (see [15])
which is equal to the identity on the smaller one.

Given two Hilbert C*-modules over a given
C*-algebra, let us say that they are equivalent
in the sense of Cuntz if, up to isomorphism, they have the same
compactly contained subobjects. In other words a third
Hilbert C*-module should be isomorphic to a
compactly contained subobject of one  of these two Hilbert
C*-modules  if and only if 
it is isomorphic to a compactly contained subobject
of the other. Clearly this is an equivalence relation, and it
can be extended to a pre-order relation by
requiring every compactly contained subobject of
the first of two given Hilbert C*-modules to be isomorphic to a
compactly contained subobject of the second. Let us
pass to the space of
Cuntz equivalence classes, so that the pre-order
becomes an order, and let us denote the resulting order
relation by $\le$.

The first comment that should be made concerning
the relation $\le$ on Cuntz equivalence classes is that
while it is  trivially an order relation (as it just refers
to comparision by inclusion of the sets of isomorphism classes of
compactly contained subobjects of two given objects),
it is not quite obvious that this relation even
holds---as it of course should!---when one object
is contained in another as a subobject. Let us
verify this. If $X$ and $Y$ are countably generated
Hilbert C*-modules with $X\subseteq Y$, and if
$X_0$ is a compactly contained subobject of $X$---i.e., 
if $X_0\subset\subset X$---, in other words, if $X_0\subseteq X$ and 
some compact self-adjoint endomorphism of 
$X$ is the identity on $X_0$, then this
compact endomorphism extends in a natural way
to a compact endomorphism of $Y\supseteq X$
(as follows from Theorem 2 of [14], which shows that
the obvious extension of a finite-rank endomorphism
has the same norm), 
which is
still the identity on $X_0$, and still self-adjoint, and so $X_0 \subset\subset Y$;
this proves that $[X]\leq [Y]$, where $[X]$ and $[Y]$
denote the Cuntz equivalence classes of $X$ and $Y$.

The second comment that should be made concerning the relation $\leq$ between 
equivalence classes is that it is compatible with addition, but that, even in 
the special case of equality, this must be proved.  In other words, even the 
fact that the Cuntz semigroup exists, in the present (nonapproximative)
context of countably generated Hilbert C*-modules, isomorphism, and compact 
containment, must form part of our analysis of this invariant. (Of course,
one approach would be just to reconcile the present equivalence relation with 
the approximative one of Cuntz---which, as is easily seen, does amount to a relation
on (countably generated) Hilbert C*-modules, weaker in general
than isomorphism. Cuntz, as it happens, only looked at
finitely generated Hilbert C*-modules,
rather than arbitrary countably generated Hilbert C*-modules,
but this is immaterial. Indeed, in the case of a stable C*-algebra,
every countably generated Hilbert C*-modules is singly generated, and
in fact is a closed right ideal, i.e., a subobject of the C*-algebra
considered as a Hilbert C*-module.
In
the general case one has to look at (countably generated)
subobjects of the countably infinite direct sum of copies of
the given C*-algebra. Note that only when the given C*-algebra
has a countable approximate unit is the Hilbert
C*-module arising from it countably generated. That
Cuntz's equivalence relation viewed as a relation
on (countably generated) Hilbert C*-modules coincides
with ours---see Appendix 6---follows from a
functional calculus lemma of Kirchberg 
and R\o rdam---Lemma 2.2 of [17]---which 
is in fact the main technical tool in our analysis below.)
\smallskip

{\bf Theorem.}\
{\it The Cuntz invariant (as defined above) is 
an ordered semigroup with zero (the order relation
compatible with addition in the sense that  if 
$[X_1]\leq [Y_1]$ and $[X_2]\leq [Y_2]$ then
$[X_1\oplus X_2]\leq [Y_1\oplus Y_2]$,
and, also, zero is the smallest element).
The order
relation has the following two purely order-theoretic properties:

(i)\ every increasing sequence---equivalently, every countable
upward directed set---has a supremum;

(ii)\  for any element the set of elements compactly
contained in it in the order-theoretic sense---where
we say that $x$ is compactly contained in
$y$ in the order-theoretic sense---let us denote this
relation by $x<<y$---if, whenever $y_1\le y_2\le \cdots$ is an
increasing sequence with supremum greater than or
equal to $y$, eventually $x\le y_n$---is upward directed---also
with respect to the stronger relation $<<$ (which
is transitive, and even antisymmetric, but not in general 
reflexive)---and contains an increasing sequence---which may be
chosen to be rapidly increasing, i.e.~with each term
compactly contained in the next---with supremum
the given element.

The operation of passing to the supremum
of a countable upward directed set and the relation
$<<$ of compact containment in the order-theoretic
sense are compatible with addition, in the sense
that the supremum of the sum $S_1+S_2$ of two
countable upward directed subsets $S_1$ and $S_2$ (which
is of course also upward directed, and so by (i)
has a supremum) is the sum of the suprema
(in other words $\sup (S_1+S_2)=\sup S_1+\sup S_2)$, and
the relations $x_1<<y_1$ and $x_2<<y_2$ imply 
$x_1+x_2<< y_1+y_2$.}
\smallskip

{\it{Proof.}}\ We must first show, in order to obtain the
Cuntz invariant as an ordered semigroup in our approach
(even just as a semigroup), that the pre-order relation
defined above on countably generated Hilbert
C*-modules (over a given C*-algebra) is compatible
with addition. Let $X_1,\, X_2,\, Y_1$, and $Y_2$ be countably
generated Hilbert C*-modules, suppose that $[X_1]\le [Y_1]$
and $[X_2]\le [Y_2]$, i.e., that every compactly
contained subobject of $X_1$ (where by subobject we
mean countably generated closed submodule, considered
with the inherited algebra-valued inner product) is
isomorphic to a compactly contained subobject of $Y_1$
and similarly for $X_2$ and $Y_2$, and let us show
that $[X_1 \oplus X_2]\le [Y_1\oplus Y_2]$. Let $X$ be a
compactly contained subobject of $X_1\oplus X_2$---i.e.,
$X\subset\subset X_1\oplus X_2$---and let us show
that $X$ is isomorphic to $Y\subset\subset Y_1\oplus Y_2$.

Consider first the case that $X=X'_1\oplus X'_2$ with
$X'_1 \subset\subset X_1$ and $X'_2 \subset\subset X_2$.
In this case, by hypothesis, $X'_1$ and
$X'_2$ are isomorphic respectively to objects 
$Y'_1\subset\subset Y_1$ and
$Y'_2 \subset\subset Y_2$, and hence $X'_1 \oplus X'_2$
is isomorphic to the subobject $Y'_1\oplus Y'_2$ of
$Y_1\oplus Y_2$,
which is clearly compactly contained in $Y_1\oplus Y_2$.
(If $b_1$ and $b_2$ are 
compact self-adjoint endomorphisms of $Y_1$ and $Y_2$ acting as the
identity on $Y'_1$ and $Y'_2$, respectively, then $b_1\oplus b_2$
is a compact self-adjoint endomorphism of $Y_1\oplus Y_2$ acting as the
identity on $Y'_1 \oplus Y'_2$.)

The case that $X$ is isomorphic to
a subobject $X'_1 \oplus X'_2$ of $X_1\oplus X_2$ as above
(i.e., with $X'_1 \subset\subset X_1$ and $X'_2 \subset\subset X_2$)
of course reduces to the preceding case (that
$X=X'_1\oplus X'_2$). The case that $X$ is
isomorphic to a subobject of such a subobject follows 
immediately. We shall show below, after first
establishing the assertions (i) and (ii) of the theorem
without making any reference to addition whatsoever,
that this last case is in fact the general case. In
anticipation of this, and in deference to the fact
that the Cuntz semigroup exists in the literature already---with a
different definition that will later be seen to be the same---let 
us in the meantime anyway refer to the ordered set of
Cuntz equivalence classes as the Cuntz semigroup.

Let us proceed to the proof of the first purely order-theoretic
assertion of the theorem (the assertion (i)).
Let $x_1\le x_2\le \dots$ be an increasing sequence
in the Cuntz semigroup of a given C*-algebra $A$,
and let us show that the supremum of the
set $\{x_n;\ n=1, 2,\cdots\}$ exists in this ordered set.
The proof consists of two steps, of which the 
first can be dealt with immediately, whereas
the second requires first establishing the second
assertion of the theorem.

The first step is to consider the case
that the sequence is rapidly increasing, i.e.,
that $x_1<<x_2<<\cdots$. In fact, this is
disingenuous, and before continuing with the
proof, we must introduce a concrete
analogue of the relation $<<$, arising from the
relation $\subset\subset$ of compact inclusion of a subobject of a
Hilbert C*-module. Namely, let us say that
$x_1$ is compactly contained in $x_2$ in the concrete sense,
and write $x_1 \subset\subset x_2$---we shall eventually show that
this relation is equivalent to $x_1 << x_2$!---if 
$x_1\le [X'_2]$ for some object (i.e., countably generated
Hilbert C*-module) 
$X'_2$ such that
$X'_2\subset\subset X_2$  
for some object $X_2$ representing $x_2$, i.e., with
$[X_2]=x_2$.
(We do not require that $x_1=[X'_2]$ for some object
$X'_2$ as above, i.e., $X'_2 \subset \subset  X_2$, with
$[X_2]=x_2$, although it might be the case that $X'_2$ can always be
chosen in this way.)

More precisely, then, the first step is to consider the case that the sequence is
rapidly increasing---written $x_1\subset\subset x_2
\subset\subset\cdots$---in the concrete sense,
i.e., in the sense that
there exist countably generated
Hilbert $A$-modules $X_1,\ X_2,\ X'_2,\ X_3,\ X'_3, \cdots$
with $X_2'\subset\subset X_2$, $X_3'\subset\subset X_3$, $\cdots$ such that 
$[X_n]=x_n$ for $n=1, 2, \cdots,$ and 
$x_1\leq [X_2']$,  $x_2\leq [X_3']$, $\cdots$ where, as above, $[X]$ denotes the Cuntz
equivalence class of $X$. Recall that, as proved above,
$X\subset\subset Y$ implies $[X]\le [Y]$, and so 
this is indeed a special case
$(x_1\subset\subset x_2\subset\subset\cdots$ does imply
$x_1\le x_2\le \cdots)$.

By the definition of Cuntz equivalence, from
$X'_2\subset\subset X_2$ and $[X_2]\le [X'_3]$ it follows that
$X'_2$ is isomorphic to a compactly contained
subobject of $X'_3$. Similarly, $X'_3$ is isomorphic
to a compactly contained subobject of $X'_4$, and so
on. Thus, we have a sequence of $A$-module
maps, preserving the $A$-valued inner product
and in particular isometric,
$$X'_2\to X'_3\to \cdots,$$
for which the image of each object is a compactly
contained subobject of the next, and we assert that
the inductive limit Hilbert $A$-module $X=\lim\limits_\to X'_n$
gives rise to the supremum of the classes 
$x_1=[X_1],\ x_2=[X_2],\ \cdots$ in the Cuntz semigroup of $A$---which are,
after all, intertwined with the classes $[X'_1]$, $[X'_2]$, $\cdots$.
In fact, we shall not need to use that the image
of $X'_n$ is a compactly contained
subobject of $X'_{n+1}$, but just a subobject for each
$n=2, 3, \cdots$.

In other words, changing notation, we must show that if $X=\lim\limits_\to X_n$ for a
sequence $X_1\rightarrow X_2\rightarrow$ of $A$-module maps preserving the
$A$-valued inner product, then $[X]=\sup [X_n]$. (As mentioned implicitly
this is what is pertinent above, since, in the earlier notation,
$[X_1]\leq [X'_2]\leq [X_2] \leq [X'_3]\leq [X_3]\leq \cdots$.) 

To show that the class of $X$ is the supremum
of $[X_1], [X_2], \cdots$, we must show that if $Y$ is a
countably generated Hilbert $A$-module such that
$[X_n]\le [Y]$ for
all $n=1, 2, \cdots$ then also $[X]\le [Y]$. Let $Y$ then\
be such that $[X_n]\le [Y]$ for all $n$ in the Cuntz
semigroup of $A$, and let $Z$ be a countably
generated Hilbert $A$-module such that $Z\subset\subset X$; we
must show that $Z\cong Z'$ for some subobject $Z'\subset\subset Y$.
We shall show that
in fact $Z\cong Z''$ for some $Z''\subset\subset X_n$ for
some $n$. (It follows 
from this by definition of the inequality $[X_n]\le [Y]$
that $Z''\cong Z'$ with $Z'\subset\subset Y$ as required.)

Let us show then that $Z\cong Z''$ for some $Z''\subset\subset X_n$.
Using the hypothesis $Z\subset\subset X$, choose
a compact self-adjoint endomorphism
$b$ of $X$ such that $b$ is
equal to the identity on $Z\subset\subset X$. We may
replace $b$ by $b^2$ to ensure that $b$ is
positive, and then replace $b$ by a
function of $b$ so that, for some $\epsilon >0$, also 
$(b-\epsilon)_+$
is the identity on $Z$. Again
replacing $b$ by a function of $b$
(in the C*-algebra generated by $b$, and so still a
compact endomorphism of $X$), we may suppose
that also yet another positive function of $b$, say $c$, is the
identity on $b$ (and still belongs to the algebra of
compact endomorphisms of X). Choose a sequence
$(c_n)$ of positive compact endomorphisms of $X$,
with $c_n$ arising from a compact endomorphism
of $X_n\subseteq X$, such that $c_n$ converges to $c$. Then
$c_n b c_n$ converges to $cbc=b$, and so, with $\epsilon >0$
as above, by Lemma 2.2 of [17], 
if $n$ is large enough that $\|c_n b c_n-b\|<\epsilon$,
then for some compact endomorphism $d_n$ of $X$
(which could be chosen to have norm at most one),
$d_n c_n b c_n d_n^*=(b-\epsilon)_+$. Note
that $c_nbc_n$ also
arises from a compact endomorphism of $X_n$ (as
these constitute a hereditary sub-C*-algebra of
the compact endomorphisms of $X$), and so also
its positive square root, $g_n$, does. Since
$$(d_n g_n)(d_n g_n)^*=(b-\epsilon)_+,$$
the partially isometric part of $d_n g_n$ in the bidual
of the C*-algebra of compact endomorphisms of $X$
therefore determines an isomorphism between a
subobject of $X_n$ and the subobject of $X$
generated by $(b-\epsilon)_+ X$ (i.e., the closure of this
submodule of $X$---recall that $A$ acts on $X$ on the right). Since
$(b-\epsilon)_+$ acts as the identity on $Z$, the object
$Z$ is a subobject of the subobject of $X$
in question, the closure of the range of $(b-\epsilon)_+$, and
is therefore isomorphic, as desired, to a subobject
of (the above subobject of) $X_n$.

This completes the first step
of the proof of the assertion (i), namely, the
consideration of the special case of a rapidly
increasing sequence $x_1 \subset\subset x_2 \subset\subset \cdots$ 
in the concrete sense---cf.~above---namely---as it turns out, but this
step was slightly subtle---, that
$x_n=[X_n]$ with $X_1\subset\subset X_2\subset\subset\cdots$. As
mentioned above,
before proceeding to the second step (the general case!)
we must first establish the second assertion of the
theorem. Actually, this is not quite correct; it
will suffice to establish the assertion (ii) with
the purely order-theoretic relation $x<<y$ replaced by the concrete
relation $x\subset\subset y$
introduced above $(x\le [Y']$ for some $Y'\subset\subset Y$ with
$[Y]=y)$.
This will allow us to complete the proof of
the assertion (i), and it will 
then be possible to deduce that the two relations $\subset\subset$ 
and $<<$---the concrete and the abstract---are equivalent.

Let, then, a countably generated (right) Hilbert
C*-module $X$ over $A$ be given. Choose, as follows, an
increasing sequence
$$X_1\subset\subset X_2\subset\subset \cdots \subset\subset X$$
of subobjects,  rapidly increasing
in the concrete sense---that each subobject $X_n$ in the
sequence is compactly contained in the next in what might be called the strong
concrete sense---note
that as pointed out above it follows that
each $X_n$ is compactly contained in $X$---such that
$X$ is the subobject generated by $X_1\cup X_2\cup\cdots$,
in other words such that the union is dense.
In fact, that $X$ is countably generated is
equivalent to the condition on the C*-algebra of compact
endomorphisms of $X$ that it have a countable
approximate unit
(see Corollary 1.1.25 of [12]).
(If ($\xi_i)$ is a generating sequence for
$X$, with $\|\xi_i\|=2^{-i}$, then the endomorphism
$\sum \xi_i\xi^*_i$ where $\xi^*_i$ denotes the module map
$\xi\mapsto\langle \xi_i,, \xi\rangle$ from
$X$ to $A$---and $\xi_i \xi^*_i$ the map $\xi\mapsto \xi_i
\langle\xi_i, \xi\rangle$!---is 
a strictly positive element of the C*-algebra of
compact endomorphisms of $X$: If $f$ is a positive functional
zero on $\xi_i\xi_i^*$ for all $i$ then, since $\xi_i a a^*\xi_i^*\le
\|a\|^2\xi_i \xi^*_i$, the corresponding complex-valued
inner product on $X$ is zero on $\xi_i A$ for every $i$ and therefore
zero on all of $X$---so $f$ is zero on every $\xi\xi^*$ and 
therefore zero. The converse is not needed and so we
leave it to the reader; it can be proved by viewing the module as
a subobject of the direct sum of a suitable number of
copies of $A$.)
Choose a countable approximate unit $(u_n)$ for the
C*-algebra of compact endomorphisms of $X$ such that 
$u_{n+1} u_n=u_n$ for every $n$. Then the increasing sequence
of subobjects
$$X_1=(u_1 X)^{-} \subseteq X_2=(u_2 X)^{-}
\subseteq \cdots \subseteq X$$
has the required properties $(X_n \subset\subset X_{x+1}$ and
$\bigcup X_n$ dense in $X$).

By the case of the assertion (i) established above,
$[X]=\sup [X_n]$ in the Cuntz (ordered)
semigroup. In fact, in the proof of this an (apparently)
stronger statement was obtained: if
$Z\subset\subset X$ then $Z$ is isomorphic to a subobject
of $X_n$ for some $n$, and hence for all sufficiently large $n$.
It follows immediately that any two compactly contained
subobjects of $X$ 
are isomorphic to subobjects
of another compactly contained subobject of $X$,
and (recall that $X_n\subset\subset X_{n+1}$ so that
$Y\subseteq X_n$ implies $Y\subset\subset X_{n+1}$)
in fact to compactly contained subobjects of such a subobject,
and hence that in the Cuntz semigroup the
set of elements $\subset\subset [X]$
is upward directed with respect to the
relation $\subset\subset$. At the same time one sees that the
supremum of this upward directed set of elements is
$[X]$. (Note that this particular upward directed
set of elements of the Cuntz semigroup contains a
cofinal increasing sequence; in general, this
is not the case, as can be seen by considering
an uncountable  direct sum of C*-algebras. One sees at the same
time that an arbitrary upward
directed subset of the Cuntz semigroup may not have a supremum.)

Now let us return to the general case of the first assertion of 
the theorem, i.e., that the sequence
$$x_1\leq x_2\leq \cdots$$
is an arbitrary increasing sequence in the Cuntz semigroup
of $A$. Choose objects (countably generated 
Hilbert C*-modules) $X_1, X_2, \cdots $ representing the Cuntz
equivalence classes $x_1, x_2, \cdots$ and as above for each
$n=1, 2, \cdots$ choose an increasing sequence of subobjects
of $X_n$, each compactly contained in the next,
$$X_{n1}\subset\subset X_{n2} \subset\subset \cdots \subset\subset X_n,$$
such that the sequence $X_{n1}, X_{n2}, \cdots $ generates $X_n$,
so that as shown above
$$[X_n]=\sup_m [X_{nm}].$$
Now note that,
as another instance of what was proved above,
concerning the supremum of an increasing sequence
of objects each contained compactly as a subobject of
the next, for each $n$ and each $m$, there exists $r$ such
that $X_{nm}$ is isomorphic to a compact subobject of
$X_{n+1, r}$. We may therefore choose $m_n$ for each $n$
in such a way that
$$\eqalign{
& [X_{11}]\subset\subset [X_{2, m_1}],\cr
& [X_{12}], [X_{2, m_1}] \subset\subset [X_{3, m_2}],\cr
& [X_{13}], [X_{2, m_1+1}], [X_{3, m_2}]\subset\subset [X_{4, m_3}],\cr}
$$
and, in general, for each $n$, the equivalence class
$[X_{n+1, m_n}]$ compactly contains (in the concrete sense) 
$[X_{n, m_{n-1}}]$ and also one new term in each of the
preceding sequences. In this way one obtains a rapidly increasing sequence
of equivalence classes,
$$[X_{1,1}]\subset\subset[X_{2, m_1}]\subset\subset
[X_{3, m_2}]\subset\subset\cdots,$$
which eventually is greater than or equal  to  each
fixed equivalence class $[X_{mn}]$, 
and each term of which
is less than or equal to some term of the
sequence $x_1\le x_2\le \cdots$. The supremum of this
sequence, which exists by the
special case of the assertion (i)
established earlier, is therefore also the supremum
of the sequence $x_1\le x_2\le \cdots$. (For each fixed $n$ it
is greater than or equal to $x_n$ because it is greater than or 
equal to
$[X_{nm}]$ for every $m$, and $x_n=\sup_m [X_{nm}]$.)
This completes the proof of the assertion (i).

Let us now prove (using the assertion (i)) that
the two notions of compact inclusion
in the Cuntz semigroup, one defined concretely as
the quotient relation arising from the inclusion up to
isomorphism of all compactly contained subobjects of one countably
generated Hilbert C*-module as subobjects of 
a fixed compactly contained subobject of  another,
written $\subset\subset$, and the other, written $<<$, defined
purely order-theoretically, in other words, as may be done in any ordered set,
by saying that one element is compactly contained in
another if, whenever the larger one is less than or equal to
the supremum of some increasing sequence
(assumed to exist in the given case), the smaller
one is already less than or equal to one of the terms
in the sequence.

Let us show first that if $x\subset\subset y$ holds (the concrete
relation) then $x<<y$ holds (the abstract relation). That 
$x\subset\subset y$ 
holds means that, with $y=[Y]$ (any choice of $Y$),
$x\le [Y']$ for some $Y'$ with $Y'\subset\subset Y$.
Let $[Y_1]\le [Y_2]\le \cdots $ be such that $[Y]\le \sup [Y_n]$,
and let us show, as required, that $x \le [Y_n]$
for some $n$.

Again, we shall use not so much the assertion (i)
itself as its proof. Recall that in the proof of
the existence of the supremum of an increasing sequence
$$x_1 \le x_2\le \cdots$$
of equivalence classes in the Cuntz semigroup, a
new increasing sequence
$$x'_1 \le x'_2 \le \cdots$$
of smaller elements (i.e., with $x'_1 \le x_1, x'_2 \le x_2, \cdots)$
was constructed with the same supremum and
such that the supremum could now be constructed as
the inductive limit of an increasing sequence of objects,
$$X'_1\subseteq X'_2 \subseteq \cdots,$$
with $X'_n$ representing $x'_n$ for each $n$.
Recall, furthermore, that it was proven in this special
setting---when the supremum can be described as an
inductive limit---that any compactly
contained subobject of the inductive limit is isomorphic to a
subobject of the finite-stage object $X'_n$ for some $n$.
(Note that this second
statement was in fact used in the proof of the first
statement; it is clearly of fundamental importance. It
should perhaps be pointed out again that, from
a technical point of view, this statement is
essentially just Lemma 2.2 of [17].)

Accordingly, passing to a new sequence with smaller terms but
the same supremum, instead of just considering a given sequence 
$$[Y_1]\le [Y_2]\le \cdots \le \sup [Y_n]$$
we may suppose that
$$Y_1\subseteq Y_2\subseteq \cdots$$
and $\sup[Y_n]=[\lim\limits_\to Y_n]$. Then, as shown
above (and just now recalled), if as
assumed above, $x \subset\subset [Y]$ and $[Y]\le \sup[Y_n]=
[\lim\limits_\to Y_n]$, so that in particular by definition of 
the relation $\subset\subset$ in the Cuntz semigroup, $x \leq [Y']$ for some
$Y'\subset\subset Y$, with therefore $Y'$ isomorphic to some
$Y''\subset\subset \lim\limits_\to Y_n$, then $Y''$ is isomorphic to a subobject
of $Y_n$ for some $n$, and since $x \leq [Y'']$, therefore $x\leq [Y_n]$ as desired.

Conversely, let us show that if $x<<y$ holds (the
abstract relation, defined just in terms of the
order relation $\le$) then $x\subset\subset y$ holds
(the concrete relation), i.e., $x \le [Y']$ for some $Y'\subset\subset Y$ with
$[Y]=y$.  By the assertion (ii), with the relation
$<<$ replaced by the relation $\subset\subset$---in which form
this assertion has now been proved---, as shown above
we may express a given element $y$ of the Cuntz semigroup as the equivalence class of the
inductive limit of a rapidly increasing sequence of objects, say
$$Y_1\subset\subset Y_2 \subset\subset\cdots.$$
Then, with $Y=\lim\limits_\to Y_n$, so that $[Y]=y$, if
$x<<y$ then, by definition, $x\le [Y_n]$ for
some $n$. By construction, $Y_n\subset\subset Y,$ as desired.

Before proceeding to the proof of the last assertion
of the theorem, concerning the very strong compatibility
of the order relation on Cuntz equivalence classes with
addition, we must return
to the problem of the very basic compatibility of the pre-order
relation on (isomorphism classes of ) countably generated
Hilbert C*-modules with addition, put aside at the
beginning of the proof. The problem was
reduced to the question of showing that if $X\subset\subset X_1
\oplus X_2$ then $X$ is isomorphic to a subobject of
$X'_1 \oplus X'_2$ for
some compactly contained subobjects $X'_1\subset\subset X_1$
and $X'_2 \subset\subset X_2$. 
Express $X_1$ and $X_2$, as described above, as the closure of
the union of rapidly increasing sequences of subobjects
$$X^1_1\subset\subset X^2_1 \subset\subset \cdots \subset X_1\
{\text{and}}\ X^1_2 \subset\subset X^2_2 \cdots\subset\subset X_2,$$
and note that then $X_1\oplus X_2$ is the closure of the union
of the rapidly increasing sequence of subobjects
$$X^1_1\oplus X_2^1 \subset\subset X^2_1 \oplus X^2_2 \subset\subset
\cdots\subset\subset X_1 \oplus X_2.$$
As shown above, $[X_1\oplus X_2]$ is then the supremum
of the increasing sequence\break
\noindent
$([X^n_1 \oplus X^n_2])$ in the
ordered set of Cuntz equivalence classes. As
also shown above, the relation $\subset\subset$
between modules implies the purely order-theoretic
relation $<<$ between their Cuntz equivalence classes,
and so, since $X\subset\subset X_1\oplus X_2$, we have
$[X]<<[X_1\oplus X_2]$. Hence,
$$[X]\le [X^n_1 \oplus X^n_2]$$
for some $n=1, 2, \cdots$. Since (again cf.~above)
there exists a subobject $X'$ of $X_1 \oplus X_2$ with
$$X\subset\subset X'\subset\subset X_1\oplus X_2,$$
we may conclude by replacing $X$ by $X'$ above that
$$[X']\le [X^n_1 \oplus X^n_2]$$
(for some $n$), from which it follows, by definition, that $X$ is isomorphic 
to a subobject of $X^n_1\oplus X^n_2$ (compactly
contained, but we don't need this), and so $X^n_1$ and
$X^n_2$ fulfil the requirements for $X'_1$ and $X'_2$.

At this point we could also adduce another proof of
additivity of the order relation $\le$ on Cuntz equivalence
classes, by showing that the relation $[X]\le [Y]$ as
defined in the present article
is equivalent to the approximative comparison relation
considered by Cuntz (from
which additivity would be immediate, just as for the
stronger order relation that one module is isomorphic to 
a submodule of another). (Incidentally, by Theorem 3 below this
last, extremely simple, (pre)-order relation is the same as the 
Cuntz pre-order relation on countably generated Hilbert
C*-modules in the case that
the C*-algebra has stable rank one---and, as mentioned above, 
it is antisymmetric---i.e., an
order relation---on isomorphism classes.) This approximative
order relation is most easily stated in the present context
in terms of compact homomorphisms: There should exist
a compact homomorphism from $Y$ to $X$
with image containing approximants to a given finite
subset of $X$. An alternative formulation is that there 
should exist a compact homomorphism from $X$ to $Y$
approximately preserving the norms of a given finite
set of elements of $X$. It is in
fact immediate from the fact that any finite
subset of $X$ is approximately contained in a compactly
contained subobject (shown above) that both types
of maps exist (from $Y$ to $X$ and from $X$ to $Y$),
when $[X]\le [Y]$,---and that
these may be chosen to be contractions. (Just use that
for any compactly contained subobject, of either $X$ or $Y$
(of $Y$ to obtain a map to $X$, and of $X$ to obtain a
map to $Y$), there exists a compact contraction into this subobject
which is approximately the identity on it.) That the existence
of one or the other kind of map (from $Y$ to $X$ or
from $X$ to $Y$) implies $[X]\le [Y]$ follows from
Lemma 2.2 of [17]. (The present notion of comparability
is shown to be equivalent to Cuntz's one in Appendix 6, below.)

To deal with the final part of the statement of the theorem,
let us show first that $\sup(S_1+S_2)=\sup S_1+\sup S_2$
if $S_1$ and $S_2$ are countable upward directed subsets
of the Cuntz semigroup of $A$. As we have shown,
$\sup S_1$ and $\sup S_2$ may both be represented
by the inductive limits of increasing sequences of countably
generated Hilbert C*-modules---the
equivalence classes of which are each less than
or equal to the equivalence class of some element
of $S_1$ or $S_2$, respectively. Clearly, the direct sum of these
inductive limits is the inductive limit of the direct sums,
and we have also shown that the inductive limit of
any increasing sequence of countably generated Hilbert
C*-modules gives 
rise to the supremum in the Cuntz semigroup.
(Both these statements were established in the course of the
proof of the assertion (i).) 
The statement follows (given that, in the Cuntz semigroup,
the relation $\le$ is compatible with addition---in other
words, one has an ordered semigroup).

Finally, we must show that the
relation $<<$ of compact containment in the order-theoretic
sense is compatible with addition, in other words,
that if $x_1 << y_1$ and $x_2 << y_2$, then $x_1+x_2 << y_1+y_2$. This
is seen immediately to hold with the
relation $\subset\subset$ in place of the relation $<<$ 
(given that the Cuntz semigroup is an ordered semigroup),
and we have proved above that these two 
relations are equivalent.
\medskip

{\bf 2.} Let us denote by $\cCu$ the category of ordered abelian
semigroups with the properties established in Theorem 1, with, as
maps, semigroup maps preserving the zero element, preserving
suprema of countable upward directed subsets, and preserving the relation
$<<$ of compact containment in the order-theoretic sense.
\smallskip

{\bf Theorem.} {\it The Cuntz semigroup is a functor from the 
category of C*-algebras, with *-homomorphisms
as maps, to the category $\cC u$, preserving inductive limits of
sequences---which always
exist in the category $\cC u$ (as well as in the category of C*-algebras).}

\smallskip
{\bf \it Proof.} Let us first show that sequential inductive limits 
exist in the category $\cC u$.  Let
$S_1 \rightarrow S_2 \rightarrow \dots$ be a sequence in the category $\cC u$.
In order to construct the inductive limit of this sequence let us first show
that the collection of increasing sequences ($s_1, s_2, \cdots$) with $s_1 \in S_1,
s_2 \in S_2, \cdots$---increasing in the sense that for each $i$ the image
of $s_i \in S_i$ in $S_{i + 1}$ is less than or equal to $s_{i+1}$---becomes
a pre-ordered abelian semigroup with the addition operation
\vskip -0.1truein
$$(s_i) + (t_i) = (s_i + t_i)$$

\noindent
and the pre-order relation

\noindent
$(s_i) \le (t_i)$ if for any $i$ and any $s \in S_i$ with $s << s_i$, 
eventually $s << t_j$ (in $S_j$).
\smallskip

Note that addition makes sense: if the sequences $(s_i)$ and $(t_i)$ with
$s_i, t_i \in S_i$ are increasing, i.e., if $s_i \le s_{i + 1}$ and
$t_i \le t_{i + 1}$ in $S_{i + 1}$ for every $i$, then by compatibility
of the order relation with addition in $S_{i + 1}$, also $s_i + t_i \le s_{i + 1} + t_{i + 1}$
for every $i$, so that the sequence $(s_i + t_i)$ belongs to the collection considered.  That
this collection becomes an abelian semigroup with this addition follows immediately
from the fact that each $S_i$ and hence also the Cartesian product
$\Pi S_i$ is an abelian
semigroup (in other words that addition in $S_i$ is associative and commutative).
We must check that the relation $(s_i) << (t_i)$ is a pre-order
relation, i.e., is reflexive and transitive.  It is reflexive because if $(s_i)$ is an 
increasing sequence with $s_i \in S_i$ and if $s << s_i$ for some $i$, then
$s << s_j$ for all $j \ge i$ (as both $s << s_i$ and $s_i \le s_j$ in $S_j$ for all
$j \ge i$, because the maps in the sequence $S_1 \rightarrow S_2 \rightarrow \cdots$
preserve both the relations $\le$ and $<<$, and since by definition if
$x << y$ and $y \le z$ in an ordered set then $x << z$).  It is transitive as an
immediate consequence of the definition (if increasing sequences $(s_i), \ (t_i)$, and
$(u_i)$, with $s_i, t_i, u_i \in S_i$, are given, such that $(s_i) \le (t_i)$ and
$(t_i) \le (u_i)$, then for any $i$ and any $s \in S_i$ with $s << s_i$, first,
eventually $s << t_j$, and in particular $s << t_j$ for some $j$ (this is in fact
the same thing), and, hence, second, eventually $s << u_k$, as is needed to show
$(s_i) \le (u_i)$).  Finally, we must check that this pre-order relation is
compatible with addition, i.e., that if $(s_i), (s'_i)$ and $(t_i), (t'_i)$ are 
increasing sequences with $s_i, s'_i, t_i, t'_i \in S_i$ and with $(s_i) \le (s'_i)$
and $(t_i) \le (t'_i)$, then $(s_i + t_i) \le (s'_i + t'_i)$.  Let $s \in S_i$ for some
$i$ be such that $ s << s_i + t_i$ and let us show that eventually $s << s'_j + t'_j$.
Choose increasing sequences
$$s^1_i \le s^2_i \le \dots << s_i \ \text{and} \ t^1_i \le t^2_i \le \dots << t_i$$

\noindent
in $S_i$ with suprema $s_i$ and $t_i$ respectively---these exist by hypothesis
(even rapidly increasing, but we do not need the full force of this).  Then as
by hypothesis the relation $\leq$ and the operation of passing to the supremum of an
increasing sequence are compatible with addition in $S_i$, we have
$$s^1_i + t^1_i \le s^2_i + t^2_i \le \dots \le s_i + t_i$$

\noindent
and furthermore $s_i + t_i = \sup_n (s^n_i + t^n_i)$.  Hence as $s << s_i + t_i$, for 
$n$ sufficiently large, $s \le s^n_i + t^n_i$.  Since $s^n_i << s_i$ and
$t^n_i << t^1_j$, it follows from $(s_k) \le (s'_k)$ and $(t_k) \le (t'_k)$
that for sufficiently large $j$ (with $i$ fixed as above) also $s^n_i << s'_j$
and $t^n_i << t'_j$. Hence by compatibility of the relation $<<$ with addition, 
for large $j$ 
$$s^n_i + t^n_i \le s'_j + t'_j,$$

\noindent
and so, combining this with $s \le s^n_i + t^n_i$, we have $s \le s'_j + t'_j$ as desired.  (It is 
immediate from the definition of the relation $<<$ in an ordered set that if $x \le y$ and $y << z$
then $x << z$---for that matter also that if $x << y$ and $y \le z$ then $x \le z$, and we shall
also use this, in the very next step.)  

Let us now show that the quotient of the pre-ordered abelian semigroup just defined by the
equivalence relation derived from the pre-order---i.e., $s$ equivalent to $t$ if $s \le t$ and
$t \le s$---is an ordered abelian semigroup belonging to $\cC u$, and furthermore is the 
inductive limit in this category of the sequence $S_1 \rightarrow S_2 \rightarrow \cdots$.

Denote this ordered abelian semigroup by $S$. To show that $S$ belongs to the category $\cC u$ we 
must show that $S$ has a zero element, that each increasing sequence in $S$ has a supremum 
(equivalently,
each countable upward direct set has a supremum) that each element of $S$ is the supremum of a
rapidly increasing sequence (each term compactly contained in the next in the order-theoretic
sense), and, finally, that the relations $\le$ and $<<$ and the operation of passing to the supremum
of an increasing sequence are compatible with addition.

The sequence $(0, \ 0, \cdots)$---or, rather, its equivalence class---is a zero element.  
(Necessarily unique.)
In order to establish the other desired properties of $S$ it is convenient to show 
first that every increasing sequence
$(s_1,\ s_2, \cdots)$ with $s_i \in S_i$ is equivalent to a rapidly increasing one.  This
is done by choosing for each $i$ a rapidly increasing sequence in $S_i$ with supremum $s_i$
(which exists by hypothesis), and then passing to a subsequence of each of these sequences, one
after another, starting with the second one, using the compact containment of each term in the
sequence for $s_i$ in $s_i$, to ensure that, for each $i$, each term of the sequence
for $s_{i+1}$ is greater than or equal to the corresponding term of the sequence for $s_i$.
Then the Cantor diagonal sequence---the $i$th term of which is the $i$th term of the (new)
sequence for $s_i$---is rapidly increasing and equivalent to the sequence $(s_i)$.  
The $(i+1)$st term of this sequence, which is the $(i+1)$st term of the sequence for $s_{i + 1}$,
majorizes the $(i+1)$st term of the sequence for each of $s_1, \cdots, s_i$, and, in particular,
as the $(i+1)$st term of the sequence for $s_i$ compactly contains the $i$th term of this
sequence, which is the $i$th term of the diagonal sequence, it follows that the $i$th term
of the diagonal sequence is compactly contained in the $(i+1)$st term.  Since the $i$th term of
the diagonal sequence is less than or equal to $s_i$, and (for the second, but not the last, time)
since $x << y$ and $y \le z$ implies $x << z$, if an element $s$ of $S_i$ is compactly contained
in the $i$th term of the diagonal sequence (an element of $S_i$), then it is compactly contained
in $s_i$ (in $S_i$, and hence by preservation of compact containment also in $S_j$ for $j \ge i$, and
hence it is also compactly contained in $s_j$ for $j \ge i$).  To prove that the diagonal 
sequence
is equivalent to $(s_i)$ it remains to show that if $s << s_i$ for some $s \in S_i$ then $s$ 
is compactly contained in all except finitely many terms of the diagonal sequence.  Choose $j$ 
such that
the $j$th term of the sequence chosen for $s_{i+1}$ is greater than or equal to $s$;
such $j$ exists by the definition of compact containment, as the supremum of this sequence, 
i.e., $s_{i+1}$,
is greater than or equal to $s_i$, also in $S_{i+1}$, and so compactly contains $s$.

Let us now show that each increasing sequence in $S$ has a supremum.  If $s^1 \le s^2 \le \cdots$ 
is an increasing
sequence in $S$, by what we have just shown there exist rapidly increasing sequences $(s^1_n)$,
$(s^2_n), \cdots$, with $s^i_n \in S_n$, representing $s^1, \ s^2, \cdots$.  
Passing to subsequences very much as above we may suppose that 
$$s^1_i << s^2_i << \cdots$$

\noindent
for each $i$, from which it follows, first (immediately), that the diagonal sequence is rapidly
increasing, i.e.,
$$s^1_1 << s^2_2 << \cdots \ ,$$

\noindent
with the $n$th inequality holding in $S_{n+1}$,
and, second, as we shall now show, that this sequence represents the supremum 
of $s^1, \ s^2, \cdots$ in $S$.  To see that the class $s$ of $(s^1_1, \ s^2_2, \cdots)$ 
is the supremum of $s^1, \ s^2, \cdots$
in $S$, note first that, for each $i$ and each $n$, $s^i_n \le s^k_k$ where $k = \max (i, n)$,
whence $s^i \le s$, and, second, that if $t = (t_1, \ t_2, \cdots) \ge s^1, \ s^2, \cdots,$ then 
for each $i$,
if $r << s^i_i$ in $S_i$ then, as $s^i_i$ is the $i$th term of $s^i$, eventually 
$r << t_j$---and so $s \le t$.

Let us show next that each element of $S$ is the supremum of a rapidly increasing sequence.  By what
we have shown above, a given element $s$ of $S$ is represented by a rapidly increasing
sequence $(s_i)$ with $s_i \in S_i$, i.e., a sequence with $s_1 << \ s_2 << \cdots$ (not just
$s_1 \le s_2 \le \cdots)$.  The sequence 
$$(s_1,  s_1, \cdots), \ (s_1,  s_2,  s_2, \cdots), \ 
(s_1,  s_2,  s_3,   s_3, \cdots),  \cdots$$  

\noindent
in $S$ is then rapidly increasing and has supremum $s$.
(It is rapidly increasing because $s_i << s_{i+1}$, not only in $S_{i+1}$ but also in $S_j$ for
$j \ge i+1$---note that preservation of compact containment by morphisms has already been used
above.  To see that the supremum of the sequence is equal to $s$, let $t \in S$ be given with
$$t \ge (s_1, s_1, \cdots), \, (s_1, s_2, s_2, \cdots), \cdots,$$

\noindent
and let us prove that $t \ge s$.  Choose just any (increasing) sequence $(t_i)$, with
$t_i \in S_i$, representing $t$.  For each $i$, we have $s_i << s_{i+1}$ in $S_{i+1}$ and
hence, as $(s_1, \cdots, s_i, s_{i+1}, s_{i+1}, \cdots) \le t$, $s_i << t_j$ in $S_j$ for all
sufficiently large $j$.  This shows that $(s_1, s_2, s_3, \cdots) \le (t_1, t_2, t_3, \cdots)$, and
so $s \le t$ in the quotient ordered set $S$.

Next, let us show that the relations $\le$ and $<<$ and the operation of passing to the supremum of an
increasing sequence in $S$ are compatible with addition.  First, recall what we have shown above, that
any single element of $S$ can be represented by a rapidly increasing sequence $(s_1, s_2, \cdots)$ with
$s_i \in S_i$, and that, furthermore, for any increasing sequence $(s^i)$ in $S$ with supremum $s$ there
is a rapidly increasing sequence $(s_1, s_2, \cdots)$, with $s_i \in S_i$, representing $s$, such that
$s_i \le s^i$ for every $i$, so that, in particular, $\sup s_i = \sup s^i$; note that it follows
from the definition of the order relation on $S$ that $\sup s_i = s$ for any rapidly increasing
sequence $(s_1, s_2, \cdots)$ with $s_i \in S_i$ representing $s$ (indeed, the construction of such a
representing sequence above shows that this is true even for a representing increasing sequence
which is not rapidly increasing).  In particular, choosing such representing sequences $(s_1, \ s_2, \cdots)$
and $(t_1, \ t_2, \cdots)$ for the suprema $s$ and $t$ of two increasing sequences $(s^i)$ and $(t^i)$ in
$S$, note that $(s_1 + t_1, s_2 + t_2, \dots)$ is a representing sequence for $s+t$ with analogous
properties---rapidly increasing, and with $s_i + t_i \le s^i + t^i$---where we do not assume that
$\sup (s^i + t^i) = s+t$, but we may now compute as follows:
$$s+t = \sup (s_i + t_i) \le \sup (s^i + t^i) \le \sup s^i + \sup t^i = s+t,$$

\noindent
which proves that 
$$\sup (s^i + t^i) = \sup s^i + \sup t^i,$$

\noindent
so that we have proved that taking suprema is compatible with addition.

>From the compatibility of the operation of taking suprema of increasing sequences in $S$
with addition, the compatibility of the relation $\le$ with addition follows.  Indeed, if 
$s^1 \le s^2$
and $t^1 \le t^2$ in $S$ then, choosing rapidly increasing representative sequences
$(s^1_i), \ (s^2_i), \ (t^1_i)$, and $(t^2_i)$ for $s^1, \ s^2, \ t^1$, and $t^2$, 
with $s^1_i, \ s^2_i, \
t^1_i$ and $t^2_i \in S_i$ (actually, it is enough for the representing sequences for $s^1$ 
and $t^1$ just to be 
increasing), we may replace $(t^1_i)$ and $(t^2_i)$ by subsequences in such a way that
$s^1_i \le s^2_i$ and $t^1_i \le t^2_i$ for every $i$, and then we have
$$\eqalign{s^1 + t^1 &= \sup s^1_i + \sup t^1_i \cr
                     &= \sup (s^1_i + t^1_i) \cr
		     &\le \sup (s^2_i + t^2_i) \cr
		     &\le \sup s^2_i + \sup t^2_i \ (\text{Equality \ not \ needed \ here.})\cr
		     &= s^2 + t^2,\cr}$$

\noindent 
i.e., $s^1 + t^1 \le s^2 + t^2$, as desired.  

The compatibility of $<<$ with addition is simpler.  If $s^1 << s^2$ and $t^1 << t^2$ in $S$ then, 
choosing
rapidly increasing representing sequences $(s^2_i)$ and $(t^2_i)$, with $s^2_i, \ t^2_i \in S_i$, from 
$s^1 << s^2$ and $t^1 << t^2$ we deduce that eventually $s^1 \le s^2_i$ and $t^1 \le t^2_i$
and so, eventually,
$$s^1 + t^1 \le s^2_i + t^2_i << s^2_{i+1} + t^2_{i+1} \le s^2 + t^2,$$

\noindent
i.e., $s^1 + t^1 << s^2 + t^2$, as desired.

Now let us complete the proof of the purely order-theoretic part of the theorem by proving that
the object $S$ of the category $\cC u$ constructed above is the inductive limit in 
this category of the
given sequence $S_1 \rightarrow S_2 \rightarrow \cdots$.  
We must show that for every object $T$ in $\cC u$
and every compatible sequence of maps $S_1 \rightarrow T, \ S_2 \rightarrow T, \cdots$ 
there exists a 
unique compatible map $S \rightarrow T$:


$$\eqalign{
  & \ T\cr
\nearrow \ \ \ \ \nearrow\hskip .4truein  & \uparrow !  \cr
  S_1 \to S_2 \to \cdots \to & \ S\cr}
$$
\noindent
Of course, for this to make sense we must have maps $S_i \rightarrow S$ for all $i$ compatible 
with the maps
$S_i \rightarrow S_{i+1}$.  For each $i$, and each $s \in S_i$, note that the sequence $(0, \cdots, 
s,  s, \cdots)$,
with $0$ until the $(i-1)$st term and then $s$ from the $i$th term on, is
increasing and therefore represents an element of $S$.  For each fixed $i$ and $s \in S_i$, note 
that for
any $j \ge i$ the sequence with $0$ up to the $(j - 1)$st term and $s$ from the 
$j$th
term on is equivalent to the one with $j=i$, defined above; this follows immediately 
from the definition of
equivalence of increasing sequences with $k$th term in $S_k$ for all $k$.  
This
shows that the maps $S_1 \rightarrow S, \ S_2 \rightarrow S, \cdots$ obtained in this 
way are compatible
with the given sequence $S_1 \rightarrow S_2 \rightarrow \cdots$ (i.e., that when 
adjoined to it they
yield a commutative diagram).

The definition of a (set) map $S \rightarrow T$ is immediate if one restricts to rapidly increasing 
representative sequences for elements of $S$ (shown above always to exist).  (If $(s_1, s_2, \dots)$ 
represents $s$ with $s_i \in S_i$ and $s_1 << s_2 << \cdots,$ let us aim to map $s$ into the 
supremum
of the increasing sequence in $T$ consisting of the images of $s_1, \ s_2, \cdots$ by the maps
$S_1 \rightarrow T, \ S_2 \rightarrow T \cdots$; this of course makes sense if the sequence $(s_1, s_2, \cdots)$ 
is just increasing.  If $(s'_1, s'_2, \dots)$ is a second rapidly increasing representative 
of $s$ then, as $(s_i)$ is rapidly increasing, for each $i$ we have $s_i << s_{i+1} \le s = \sup s'_i$ and 
so eventually $s_i \le s'_j$; by symmetry, also each $s'_i$ is eventually $\le s_j$, and so the 
suprema
of the images of $s_1, \ s_2, \cdots$ and of $s'_1, \ s'_2, \cdots$ in $T$ are equal  (as the maps
$S_i \rightarrow T$ preserve the order relation and are compatible as set maps).)

Let us check that the map $S \rightarrow T$ thus defined is compatible with the maps $S_i \rightarrow T$
(i.e., that the diagram is commutative as a diagram of set maps), that it is a morphism in 
the category
$\cC u$, and that it is unique with these properties.

To show compatibility of $S_i \rightarrow T$ with $S \rightarrow T$, we must show, for each
fixed $i$, that if $s \in S_i$ then the image of $s$ in $T$ by the given map
$S_i \rightarrow T$ is the same as the image of $s$ in $T$ by the composed map
$S_i \rightarrow S \rightarrow T$.  By definition, the image of $s$ in $S$ is represented
by the sequence $(0,  \cdots,  s,  s, \cdots)$ consisting of the element $s$ repeated
beginning with the $i$th coordinate; however, in order to compute the image of this
element of $S$ in $T$ we must represent it by a rapidly increasing sequence:  let us
use the sequence $(0, \cdots,  r_i,  r_{i + 1}, \cdots)$---i.e., $r_j$ in the $j$th
place for $j \ge i$---where $(r_j)$ is a rapidly increasing sequence in $S_i$ with supremum
$s$.  By definition the corresponding element in $T$ is the supremum of the (increasing
sequence of) images of the elements $r_i \in S_i, \ r_{i + 1} \in S_{i + 1}, \cdots$ by the maps
$S_i \rightarrow T, \ S_{i+1} \rightarrow T, \cdots$, equivalently (by commutativity) of the 
images of $r_i, \ r_{i+1}, \cdots \ \in S_i$ by the map $S_i \rightarrow T$, and as this map
preserves increasing sequential suprema, the supremum in $T$ in question is just the image
of $s$, by the map $S_i \rightarrow T$, as desired.

To show that the map $S\rightarrow T$ belongs to the category $\cC u$, we must show that it
preserves addition, preserves the order relation, preserves suprema of increasing sequences, 
and preserves the order-theoretic relation
$<<$ defined earlier, in terms of the two notions just mentioned (the order relation
$\le $ and the operation of sequential increasing supremum).  Let us address these issues, 
briefly, in turn.

Given two rapidly increasing sequences $(r_1,  r_2, \cdots)$ and $(s_1,  s_2, \cdots)$ 
with $r_i, \ s_i \in S_i$ for each $i$, to check that the sum in $S$ maps into the sum of 
the images in $T$ it is enough to recall what these images are, and that the operation of
passing to the supremum of an increasing sequence in $T$ is compatible with addition
in $T$.  To check that the relation $(r_1, r_2, \cdots) \ \le \ (s_1, s_2, \cdots)$ in
$S$ (with $(r_i), \ (s_i)$ rapidly increasing as above) leads to the same relation 
between the images in $T$, recall that the equivalence relation defining $S$ is just
derived from the pre-order relation between sequences which leads to the order relation
in $S$---and we have already shown that the map $S \rightarrow T$ exists!  (And this by 
the only way conceivable, namely, by just proving that the pre-order is preserved.)

To check that the map from $S$ to $T$ preserves suprema of increasing sequences, recall from
the proof that such suprema exist in $S$, given above, that representatives of an increasing 
sequence of elements of $S$ may be chosen in such a way that, not only is each representing
sequence rapidly increasing, but also, the sequence of $i$th terms for each fixed $i$
is rapidly increasing---and then the diagonal sequence is also rapidly increasing
and furthermore represents the supremum of the given increasing sequence in $S$.  Recalling
the definition of the map $S \rightarrow T$ one sees immediately then that the image of the
supremum is the supremum of the images.  (The images of the terms of the diagonal 
sequence eventually majorize the images of the terms of each of the representing
sequences, and so their supremum majorizes the supremum of the suprema---of the images of the
terms of the individual representing sequences.)
  
Finally, to check that the map $S \rightarrow T$ is compatible with the relation $<<$,
let $(r_1, r_2, \cdots)$ and $(s_1, s_2, \cdots)$ be rapidly increasing
representing sequences as above (i.e., with $r_i,  s_i  \in  S_i)$, and suppose that
$(r_i) << (s_i)$; we must show that the supremum of the images of $r_1, \ r_2, \cdots$
in $T$ is compactly contained (of course in the order-theoretic sense) in the supremum
of the images of $s_1 \ s_2, \cdots$---in other words, that, if $\dot r$ denotes the image of
$r = \sup r_i$ in $T$, and $\dot s$ the image of $s = \sup s_i$ in $T$, then
$\dot r << \dot s$.  The proof is very simple: since $s = \sup s_i$ and $r << s$,
we have $r \le s_i$ for some $i$.  Since the map $S \rightarrow T$ is already known to
preserve the relation $\le$, it follows that $\dot r \le \dot s$ in $T$.  Since
$s_i << s_{i+1}$, not only in $S$ but also in $S_{i + 1}$, it follows from the properties of the given
map $S_{i+1} \rightarrow T$ that $\dot s_i << \dot s_{i+1}$ in $T$.  Again, since
$s_{i+1} \le s$ in $S$, we have $\dot s_{i+1} \le \dot s$ in $T$, and hence
$\dot s_i << \dot s$ in $T$, and hence also $\dot r << \dot s$ in $T$, as desired.

Finally, let us show that the association of the Cuntz ordered semigroup to a C*-algebra gives
rise in a natural way to a functor from the category of all C*-algebras to $\cC u$, and that 
this functor preserves inductive limits of sequences.

By the functor's arising naturally we just mean that we have still to describe the morphism between
the Cuntz semigroups that should correspond to a morphism between C*-algebras.  Given
C*-algebras $A$ and $B$ and C*-algebra morphism $A \rightarrow B$ (a map preserving the *-algebra
structure---necessarily a contraction), to each given Hilbert C*-module $X = X_A$ over
$A$, associate the Hilbert C*-module over $B$ defined by completing the (right) $B$-module
$(X_A) \otimes_A (_A B)$ with respect to the (possibly degenerate) $B$-valued inner product
$$\langle \sum \xi_i  \otimes  b_i,  \sum \xi'_j \otimes b'_j \rangle_B 
= \sum\limits_{i, j} b^*_i \langle \xi_i, \xi'_j \rangle_A b'_j,$$

\noindent
where $B$ is considered as a left $A$-module by virtue of the given homomorphism $A
\rightarrow B$.

Note that this correspondence, from Hilbert $A$-modules to Hilbert $B$-modules,
takes countably generated Hilbert $A$-modules to countably generated
Hilbert $B$-modules.  Let us show that it preserves (in a natural way)
the relation of inclusion as subobject, and also preserves the relation
of compact inclusion---i.e., inclusion as a compactly contained subobject.  Since
it clearly preserves the relation of isomorphism (in a natural way), it follows
that it preserves the Cuntz pre-order relation (as defined above---i.e.,
compactly contained subobjects of the first of two objects isomorphic to
compactly contained subobjects of the second).

In fact it is also clear that a morphism $X \rightarrow Y$ of Hilbert $A$-modules,
by which let us mean one preserving the $A$-valued inner product, is transformed
by the push-forward construction described above into a morphism $X_B \rightarrow Y_B$
(of Hilbert $B$-modules).  (At the purely algebraic level of the construction
it is immediate that the natural push-forward maps $(X_A) \otimes (_A B) \rightarrow
(Y_A) \otimes (_A B)$ (tensor products over $A$) preserves the $B$-valued inner
product, and hence is isometric (although at this stage the norms may be
seminorms), and then it follows by continuity that this holds for the extension
to the completion.  (One does not need this but note that by Theorem 3.5 of [15]) it 
is enough to note that the extension is isometric, and a $B$-module map, as this
by itself implies preservation of the $B$-valued inner product.)

Let, then, $X_A$ be compactly contained in $Y_A$, and let us show that the
pushed forward inclusion of $X_B$ in $Y_B$ is a compact one, i.e., that there is
a compact self-adjoint endomorphism of $Y_B$ that acts as the identity on $X_B \subseteq Y_B$.
The hypothesis is that there exists a compact self-adjoint endomorphism, say $t$, of $Y_A$ which
acts as the identity on $X_A \subseteq Y_A$. It is sufficient to show that
the endomorphism $t = t_A$ has a push-forward, in the natural sense, to a
compact self-adjoint endomorphism $t_B$ of $Y_B$ (putting $X$ aside completely).  The
natural property that $t_B$ should have is of course that 
$t_B (\eta \otimes b) = (t_A\eta) \otimes b$, for $\eta \in Y_A$ and $b \in B$.
This condition determines purely algebraically a map on $(Y_A) \otimes_A (_A B)$,
which is bounded because $\langle t_A \eta, t_A \eta \rangle_A \le \|t^*_A t_A \|
\langle \eta,  \eta \rangle_A$ (and the fact that this also holds when
$\eta$ is replaced by $\eta_1 \otimes \dots \otimes \eta_k \in Y_A \otimes
\dots \otimes Y_A$).  This shows that any adjointable endomorphism of $Y_A$
can be pushed forward to $Y_B$; that the push-forward, $t_B$, of $t_A$ 
is compact if $t_A$ is follows from the fact that this is clear (purely
algebraically) if $t_A$ is of finite rank (i.e., a finite sum of
endomorphisms $\eta \mapsto \zeta \langle \zeta', \eta \rangle$ with
$\zeta, \zeta' \in Y_A$), together with the fact that, by definition,
$t_A$ is a limit in norm of endomorphisms of finite rank (and, also, the 
fact that the push-forward of an arbitrary adjointable endomorphism is seen
by the calculation outlined above to have at most the same norm). 

This shows that the correspondence $X_A \rightarrow X_B$ is functorial, for 
a fixed map $A \rightarrow B$, in a way that passes naturally to a morphism 
$\cC u (A) \rightarrow \cC u (B)$.

Let us now show that the resulting functor, from the category of all C*-algebras
to the category $\cC u$---for it is manifestly a functor (i.e., respects
composition of maps)---, preserves (sequential) inductive limits.

Let $A_1 \rightarrow A_2 \rightarrow \cdots$ be a sequence of C*-algebras,
with inductive limit $A$.  Let us show that if $X$ is a countably 
generated Hilbert C*-module over $A$, then the class $[X]$ of $X$
in the Cuntz semigroup of $A$ is the supremum of the increasing 
sequence consisting of the canonical images in this semigroup of a sequence
$(x_i)$ with $x_i$ in the Cuntz semigroup of $A_i$ and with
$x_i \le x_{i+1}$ for each $i$, where $x_i$ denotes the image
of $x_i$ in the Cuntz semigroup of $A_{i+1}$ and the 
comparison is  in that ordered semigroup.  As we shall show
below, this makes it possible to deduce, just from the construction
of the inductive limit of a sequence in the category $\cC u$ 
given above at the beginning of the proof, that the Cuntz
semigroup of $A$---let us denote this by $\cC u (A)$---is the 
inductive limit (in $\cC u$) of the sequence $\cC u (A_1)
\rightarrow \cC u (A_2) \rightarrow \cdots$, corresponding to the
given sequence $A_1 \rightarrow A_2 \rightarrow \cdots$.

By Theorem 2 of [14], $X$ is isomorphic to a subobject of the countably 
infinite Hilbert C*-module direct sum, $\bigoplus_1^{\infty} A$, of copies of 
$A$, and so we may suppose that it
is a subobject of $\bigoplus_1^{\infty} A$.  Although $A$ itself may not be
countably generated as a Hilbert $A$-module, and so also not $\bigoplus_1^{\infty} A$,
there exists a countably generated closed submodule $A'$ of $A$
such that $\bigoplus_1^{\infty} A'$ contains $X$---for instance, that
generated by all the coordinates in $\bigoplus_1^{\infty} A$ of a countable
generating set for $X \subseteq \bigoplus_1^{\infty} A$.  Denote
$\bigoplus_1^{\infty} A'$ by $Y$.

Concerning the object $Y$, we shall use only that $Y$ contains $X$, and
that there exists a sequence $Y_1 \subseteq Y_2 \subseteq \dots \subseteq Y$
of subobjects of $Y$ such that each $Y_n$ arises from some finite stage of
the sequence $A_1 \rightarrow A_2 \rightarrow \dots \rightarrow A$, by
the functorial push-forward construction described above.  Passing
to a subsequence of $A_1 \rightarrow A_2 \rightarrow \cdots$ we may
suppose that $Y_n$ arises form the $n$th stage, say from the object 
$(Y_n)_{A_n}$ over $A_n$.  (Let us then sometimes write $Y_n$ for
$(Y_n)_{A_n}$, and $(Y_n)_A$ for the push-forward!) (To obtain the desired
increasing sequence $Y_1 \subseteq Y_2 \subseteq \dots$ of subobjects
of $Y$ over $A$, use that $Y = \bigoplus_1^{\infty} A'$, and, if 
the closed right $A$-module $A'$ is not already the closure of the union
of the push-forwards $(A' \cap A_1)A, (A' \cap A_2)A, \cdots$ of $A_1, A_2, \cdots$
respectively, then simply enlarge $A'$ (by adjoining countably many 
elements, of $A_1, A_2, \cdots$ approximating the generators of $A'$), so 
that this is the case.)  It follows in particular that every
compact endomorphism of $Y$ is the limit of a sequence of compact
endomorphisms arising from $(Y_n)_{A_n}$, for $= 1, 2, \cdots$.

Let us now construct as promised an increasing sequence $x_1 \le x_2 \le \cdots$
in $\cC u (A)$, with $x_i$ arising from an element of $\cC u (A_i)$
for each $i$, with $x_i\leq x_{i+1}$ in $\cC u(A_{i+1})$ for each $i$,
and with sup $x_i = [X]$.  Choose a strictly positive element, $h$,
of the compact endomorphism algebra of $X$, and recall that then $h$ belongs
in a natural way to the compact endomorphism algebra of $Y$, of which $X$
is assumed to be a subobject.  Recall that $Y = (\bigcup (Y_i)_A)^-$ where
$Y_i$ is a Hilbert $A_i$-module for each $i$ with $(Y_i)_{A_{i+1}} \subseteq
Y_{i+1}$, so that in particular $(Y_1)_A \subseteq (Y_2)_A \subseteq \cdots$,
and that, furthermore, every compact endomorphism of $Y$ can be 
approximated arbitrarily closely in norm by the push-forwards of compact
endomorphisms of $Y_1, Y_2, \cdots$.  Choose a compact endomorphism 
$h_n \ge 0$ of $Y_n$ such that $h_n = (h_n)_A$ tends to $h$ in the compact
endomorphism algebra of $Y$.  By Lemma 2.2 of [17], for given
$\epsilon > 0$, choosing $h_n$ close enough to $h$---strictly within $\epsilon$
is sufficient (at this stage!---see below)---we obtain
$$(h_n - \epsilon)_+ = d  h  d^*$$

\noindent
for some compact endomorphism $d$ of $Y$ (which may be chosen to have norm
one, but we don't need this).  Furthermore, inspection of the construction
of $d$ in the proof of Lemma 2.2 of [17] shows that, with $d$ constructed
in this way, the element $h^{{1 \over 2}} d^* d h^{{1 \over 2}}$, a compact
endomorphism of $X$, is close (in norm) to $h$---to obtain this it is no
longer sufficient for $h_n$ just to be within $\epsilon$ of $h$, but for a 
given desired degree of approximation by 
$h^{{1 \over 2}} d^* d h^{{1 \over 2}}$---let us say $\epsilon !$---we may just
choose $h_n$ to give the necessary approximation to $h$.  Then, again by
Lemma 2.2 of [17], 
$$(h - \epsilon)_+ = e h^{{1 \over 2}} d^* d h^{{1 \over 2}} e^*$$

\noindent
for some compact endomorphism $e$ of $X$. Combining these two 
equations, we see that $((h - \epsilon)_+ X)^-$, a (compactly contained)
subobject of $X$, is isomorphic to a subobject of $((h_n - \epsilon)_+ Y)^-$.
(The partially isometric part of the compact homomorphism
$d h^{{1 \over 2}} e^*$ is a (not necessarily adjointable) isometry
from $((h - \epsilon)_{+} X)^{-}$ to $((h_n - \epsilon)_{+} Y)^-$, with image
$(d h^{{1 \over 2}} e^* e h^{1 \over 2} d^* Y)^- \subseteq ((h_n - \epsilon)_+ Y)^-.)$
Therefore, in $\cC u (A)$, 
$$[((h - \epsilon)_+ X)^-] \le [((h_n - \epsilon)_+ Y)^-].$$

The important point is that the object $((h_n - \epsilon)_+ Y)^-$ is the closure
of the union of the increasing sequence of subobjects $((h_n - \epsilon)_+ (Y_k)_A)^-$,
each one of which arises from a finite stage---namely, $((h_n - \epsilon)_+ (Y_k)_A)^-$
is the push-forward of $((h_n - \epsilon)_+ (Y_k)_{A_l})^-$, where
$l = \max (k, n)$.  As shown earlier, the first property implies that,
in $\cC u (A)$,
$$[((h_n - \epsilon)_+ Y)^-] = \sup\limits_{k} [((h_n - \epsilon)_+ (Y_k)_A)^-].$$

At the same time, on choosing a sequence $\epsilon_m$ tending (strictly)
monotonically to $0$, one has 
$$X = (\bigcup ((h - \epsilon_m)_+ X)^-)^-,$$

\noindent
and hence, in $\cC u (A),$
$$[X] = \sup_m [((h - \epsilon_m)_+ X)^-].$$

Finally, note that, also, for each $\epsilon$, with $h_n$ as above,
the object $((h_n - \epsilon)_+ Y)^-$ is isomorphic to a subobject of $X$,
by means of the partially isometric part of the compact
homomorphism $h^{{1 \over 2}} d^*$, from $((h_n - \epsilon)_+ Y)^-$ to
$X = (h Y)^-$.  Making $\epsilon$ smaller, we may ensure that the image
is compactly contained in $X$.  Hence, as shown earlier,
$$[((h_n - \epsilon)_+ Y)^-] << [X]$$ 

\noindent
in $\cC u (A)$.  The conclusion now follows, with $x_1, x_2, \cdots$
chosen after passing to the subsequence $A_{k_1} \rightarrow A_{k_2}
\rightarrow \cdots$ of $A_1\rightarrow A_2\rightarrow\cdots$ to be 
the sequence (with respect to the given
sequence $A_1 \rightarrow A_2 \rightarrow \cdots$)
$$y'_m = [((h_{n_m} - \epsilon_{l_m})_+ (Y_{k_m})_A)^-], \ m=1, 2, \cdots$$

\noindent
(arising as observed above from finite stages), for suitable sequences
$(n_m), (l_m)$ and $(k_m)$.  Indeed, first choose $n_1$ such that
$$[((h - \epsilon_2)_+ X)^-] \le [((h_{n_1} - \epsilon_2)_+ Y)^-].$$

\noindent
Since $\epsilon_2 < \epsilon_1$, we have $((h_{n_1} - \epsilon_2)_+ Y)^- \subset\subset 
((h_{n_1} - \epsilon_1)_+ Y)^-$, whence by the  order-theoretic compact
inclusion of the corresponding Cuntz classes in $\cC u (A)$, we may
choose $k_1$ such that
$$[((h_{n_1} - \epsilon_2)_+ Y)^-] \le [((h_{n_1} - \epsilon_1)_+ Y_{k_1})^-].$$

\noindent
By compact inclusion of the latter class in $[X]$, we may choose $l_2$
such that
$$[((h_{n_1} - \epsilon_1)_+ Y_{k_1})^-] \le [((h - \epsilon_{l_2+1})_+ X)^-].$$

\noindent
Choose $n_2$ in the same way as $n_1$ above such that

$$[((h - \epsilon_{l_2+1})_+ X)^-] \le [((h_{n_2} - \epsilon_{l_2+1})_+ Y)^-].$$

\noindent
Again by compactness, as $\epsilon_{l_2+1} < \epsilon_{l_2}$, we may choose $k_2$
in the same way as $k_1$ above such that
$$[((h_{n_2} - \epsilon_{l_2+1})_+ Y)^-] \le [((h_{n_2} - \epsilon_{l_2})_+ Y_{k_2})^-].$$

\noindent
Continuing in the way just described, we obtain an increasing sequence
$y'_m = [((h_{n_m} - \epsilon_{l_m})_+ (Y_{k_m})_A)^-]$ in $\cC u (A)$
which is intertwined with respect to the order relation with the increasing
sequence $x'_m = [((h - \epsilon_{l_m+1})_+ X)^-]$---and is in particular
increasing!
Since the second sequence, $(x'_m)$, as shown earlier, has supremum
$[X]$ in $\cC u (A)$, the first sequence, $(y'_m)$, also does.  It
remains to note that after passing to the subsequence $A_{k_1} \rightarrow
A_{k_2} \rightarrow \cdots$ and changing notation, $y'_m$ belongs to 
$\cC u (A_m)$ for each $m$, and so the sequence $(y'_m)$ fulfils the 
requirements for the sequence $(x_m)$---except possibly for the condition
$y_m'\leq y'_{m+1}$ in $\cC u(A_{m+1})$. This is ensured by again applying
Lemma 2.2 of [17], to refine the choice of the sequences $(n_m),\, (k_m),$ and
$(l_m)$ so that, for each $m$,
$$(h_{n_m} - \epsilon_{l_{m}})_{+} = d(h_{n_{m+1}} -\epsilon_{l_{m+1}})_{+}d^*$$
for some compact endomorphism $d$ of $Y_{n_{m+1}}\subseteq Y_{k_{m+1}}$, from 
which it follows that $((h_{n_m}-\epsilon_{l_{m}})_+Y_{k_m})^-$ is isomorphic
to a subobject of $((h_{n_{m+1}}-\epsilon_{l_{m+1}})_+ Y_{k_{m+1}})^-$, over 
$A_{k_{m+1}}$---then one has $y'_m\leq y'_{m+1}$ in $\cC u(A_{k_{m+1}})$, 
or in $\cC u(A_{m+1})$ after
the prescribed passage to the subsequence $A_{k_1}\rightarrow A_{k_2}\rightarrow\cdots$
of $A_1\rightarrow A_2\rightarrow\cdots$, as desired.

Let us now show that, with respect to the canonical sequence in the
category $\cC u$ corresponding to the sequence of C*-algebras
$A_1 \rightarrow A_2 \rightarrow \cdots \rightarrow A$,
$$\cC u (A_1) \rightarrow \cC u (A_2) \rightarrow \dots \rightarrow
\cC u (A),$$

\noindent
$\cC u (A)$ is the inductive limit.  By the construction of 
${\lim\limits_\to} \cC u (A_i)$ in the proof of Theorem 2,
and what has just been proved, it is sufficient to show that
if $x_1 \le x_2 \le \cdots$ with $x_i \in \cC u (A_i)$ and also
$y_1 \le y_2 \le \cdots$ with $y_i \in \cC u (A_i)$, then $\sup
x_i \le \sup y_i$ in $\cC u (A)$ (the suprema of course
referring to the images of the sequences $(x_i)$ and $(y_i)$
in $\cC u (A)$) if, and only if, whenever $z << x_i$ in $\cC u (A_i)$
for some $i$ then $z << y_j$ in $\cC u (A_j)$ for some
$j \ge i$.  (This establishes an isomorphism of ordered semigroups
from ${\lim\limits_\to} \cC u (A_i)$ onto a sub
ordered semigroup of $\cC u (A)$, and this subsemigroup was shown
above to be all of $\cC u (A)$.)

Replacing $(x_i)$ and $(y_i)$ with equivalent sequences (which does 
not change the statement of what is to be proved), as in the 
proof of Theorem 1, we may suppose, first, that $(x_i)$ and $(y_i)$
are rapidly increasing, and, second, that $x_i = [X_i]$ and
$y_i = [Y_i]$ where $X_i$ and $Y_i$ are Hilbert C*-modules
over $A_i$ and 
$$X_1 \subset \subset X_2 \subset \subset \cdots, \quad Y_1 \subset \subset
 Y_2 \subset \subset\cdots.$$

\noindent
Then
$$\sup x_i = [{\lim\limits_\to} (X_i)_A] \ \text{and} \
\sup y_i = [{\lim\limits_\to} (Y_i)_A],$$

\noindent
where $(X_i)_A$ denotes the push-forward of $X_i$ from $A_i$ to
$A$, discussed above.

Suppose first that $\sup x_i \le \sup y_i$ (in $\cC u (A)$), and let
$z << x_i$ in $\cC u (A_i)$ be given for some fixed $i$.  By
the concrete definition of $<<$ (see proof of Theorem 1),
$$z \le [Z'] \ \text{for some} \ Z' \subset \subset X_i \ \text{over} \ A_i.$$

\noindent
Furthermore, we may choose $Z''$ and $Z'''$ such that
$$Z' \subset \subset Z'' \subset \subset Z''' \subset \subset X_i \ \text{over} \ A_i.$$

\noindent
Then,
$$Z'_A \subset \subset Z''_A \subset \subset Z'''_A \subset \subset (X_i)_A.$$

\noindent
In particular, $[Z'''] << [X_i] \le \sup x_i \le \sup y_i$
(in $\cC u (A)$), and so by the abstract definition of $<<$,

$$[Z'''] \le [Y_j] \ \text{for  some} \ j ,\ {\text{in}} \ \cC u (A).$$

It follows, in particular (by definition) that $Z''$ is isomorphic to 
a compactly contained subobject of $Y_j$ over $A$, say $Y'_j$.  
Let $h$ be a positive element of the algebra of compact endomorphisms
of $(Z'')_{A_i}$ (to be specified later!).  This is then $x^*x$ where
$x x^*$ is a compact endomorphism of $(Y_j)_A$, for a certain
compact homomorphism $x$ from $(Z'')_A$ to $(Y_j)_A$ (namely, the
product of $(h^{{1 \over 2}})_A$ with an isomorphism from $(Z''_{A_i})_A$
to a subobject of $(Y_j)_A$).  Since $Z''$ and $Y_j$ arise from the $i$th
and $j$th stages, respectively, for some $k \ge \max (i, j)$ we may 
approximate $x$ in norm, in the algebra of compact homomorphisms from
$(Z'')_A$ to $(Y_j)_A$, by a compact homomorphism from $((Z'')_{A_i})_{A_k}$ to
$((Y_j)_{A_j})_{A_k}$, say $x'$.  Then $x'* x'$ is close to $h$ in the 
algebra of compact endomorphisms of $(Z''_{A_i})_A$, and since both
$h$ and $x'* x'$ arise from a finite stage, say for the moment
$(Z''_{A_i})_{A_k}$ for fixed $k$, if they are close over $A$ then 
they are also (almost as) close over $A_l$ for some $l \ge k$.  Hence
by Lemma 2.2 of [17], correcting $x'$ by composing with a compact
endomorphism of $((Z'')_{A_i})_{A_l}$, we may suppose that
$x'* x' = (h - \epsilon)_+$, for a given $\epsilon > 0$---if we choose $k$ large to 
begin with (and then $l$ large depending on this choice).

Now let us use that, since $Z' \subset \subset Z''$ over $A_i$,
there exists a compact self-adjoint endomorphism $h = h_{A_i}$ of $(Z'')_{A_i}$ equal
to the identity on $Z'$.  As shown earlier, we may suppose not only
that $h$ is positive, but that also in fact $(h - \epsilon)_+$ acts as
the identity on $Z' \subset Z''$.  This gives the choice of $h$ and
$\epsilon$ to be used above.  We thus obtain that $x'$ is a compact
homomorphism from $(Z'')_{A_l}$ to $(Y_j)_{A_l}$ such that
$x'* x' = (h - \epsilon)_+$, from which we deduce that the restriction
of $x'$ to $(Z')_{A_l} \subset \subset (Z'')_{A_l}$ is an isomorphism
from $(Z')_{A_l}$ to a compactly contained subobject of $(Y_j)_{A_l}$
(as $x' x'*$ is a compact endomorphism of $(Y_j)_{A_l}$ acting
as the identity on the image of $(Z')_{A_l}$ by $x'$---given
that the compact self-adjoint endomorphism $x'* x'\,\, (= (h - \epsilon)_+)$ of
$(Z'')_{A_l}$ acts as the identity on $(Z')_{A_l})$, as desired.

Suppose, conversely, that whenever $z << x_i$ in $\cC u (A_i)$
for some $i$ then $z << y_j$ in $\cC u (A_j)$ for some
$j \ge i$, and let us show that $\sup x_i \le \sup y_i$ in
$\cC u (A)$. We must show that $x_i \le \sup y_j$ in $\cC u (A)$
for every $i$.  Let us then fix $i = 1, 2, \cdots$.  By 
Theorem 1, $x_i$ is the supremum in $\cC u (A_i)$ of an 
increasing sequence $(z_n)$ in $\cC u (A_i)$
with $z_n << x_i$ in $\cC u (A_i)$ for each $n = 1, 2, \cdots$.
Then by hypothesis, for each $n$ we have $z_n << y_j$ in
$\cC u (A_j)$ for some $j \ge i$.  Then also (by functoriality)
$z_n << y_j$ in $\cC u (A)$---where now $j$ is fixed but $n$ is
arbitrary.  In fact, all we shall need is $z_n \le y_j$.
By functoriality, from $x_i = \sup z_n$ in $\cC u (A_i)$
follows $x_i = \sup z_n$ in $\cC u (A)$, and so $x_i = \sup 
z_n \le y_j$ in $\cC u (A)$, as desired.
\smallskip

{\bf 3.}\
The following result is a partial answer to
the questions concerning Cuntz equivalence raised
in Section 1.
\smallskip

{\bf Theorem.}\
{\it Let $A$ be a C*-algebra of stable
rank one. Two countably generated
Hilbert C*-modules over $A$ are equivalent in the
sense of Cuntz (described in Section 1) if, and only if,
they are isomorphic. In other words, the Cuntz semigroup in
this case is just the semigroup of isomorphism classes
of countably generated Hilbert C*-modules.
Furthermore, the order structure arises from
inclusion of modules.}
\smallskip

{\it Proof.}\
Since isomorphic Hilbert C*-modules have
the same isomorphism classes of compactly contained
subobjects, by definition they are Cuntz equivalent.

Let $X$ and $Y$ be Cuntz equivalent countably
generated Hilbert C*-modules over $A$, and let us show that
$X$ and $Y$ are isomorphic. Choose (as described above)
rapidly increasing sequences of subobjects
$$X_1\subset\subset X_1\subset\subset \cdots \subseteq X,$$
$$ Y_1 \subset\subset\, Y_2 \subset\subset \cdots \subseteq Y,$$

\noindent
generating $X$ and $Y$ respectively. By the definition of Cuntz
equivalence, $X_2$ is isomorphic to a compactly contained
subobject of $Y$. In particular (note that this is a priori
a weaker property), $[X_2]\subset\subset[Y]$ 
(i.e., $[X_2]\le [Z]$ for some $Z\subset\subset Y$),
and hence as shown in the proof of Theorem 1,
$[X_2]<<[Y]$. As shown in the proof of Theorem 2,
$[Y]=\sup [Y_i]$. It follows from the definition of compact
containment (of Cuntz equivalence classes) in the order-theoretic
sense that $[X_2]\le [Y_i]$ for some $i=1, 2,\cdots$. In
particular, on choosing $i_1$ with $[X_2]\le [Y_{i_1}]$, we have
by definition that, as $X_1\subset\subset X_2$, the object $X_1$ is
isomorphic to a subobject of $Y_{i_1}$, say by the
map $\varphi_1$: $X_1\to Y_{i_1}$.

In a similar way (considering first $Y_{i_1+1} \subset\subset Y)$
we obtain an isomorphism $\psi_1$ of $Y_{i_1}$ onto a
subject of $X_{j_1}$ for some $j_1=1, 2, \cdots$. Continuing in
this way, and passing to subsequences of $(X_i)$
and $(Y_i)$ and changing notation, we have a diagram
$$\eqalign{
& \quad X_1 \ \subset\subset\ \ X_2 \ \ \subset\subset \cdots \subseteq X\cr
& \varphi_1 \downarrow \psi_1\nearrow{\varphi_2}\downarrow
{\psi_2}\nearrow\cr
& \quad Y_1 \ \ \subset\subset \ \  \ Y_2\ \ \subset\subset\cdots \subseteq Y,\cr}
$$
in which each vertical map (downwards or
upwards) is an isomorphism onto its image.

It is sufficient, by a modification of 2.2 and
2.1 of [6] (using that $X_i$ and $Y_i$ are countably
generated)---see Example 4.4 and Theorem 3 of [7]---to show that any
two isomorphisms from one Hilbert C*-module over A onto
submodules of another one are approximately equal,
on finitely many elements,
modulo inner
automorphisms of the codomain Hilbert module,
i.e., automorphisms arising from unitary elements
of the C*-algebra of compact endomorphisms with
the identity adjoined. (In other words, any two such 
homomorphisms are close on finitely many elements after one
of them is composed with such an automorphism.)

Let us establish this fact, using, naturally, that
$A$ has stable rank one. It is enough to show that,
in this case, an isomorphism between
two closed submodules  of a Hilbert C*-module can be
approximated pointwise by 
an inner automorphism (defined as above) of the
whole module. The first step is to note that,
by Proposition 1.3 of [15], such an isomorphism
can be approximated on each finite set by a compact
homomorphism of norm one and that---cf.~above---such a homomorphism
extends to a compact endomorphism of
norm one of the larger Hilbert module. 
The second step is to
note that, as
the property of having stable rank one is invariant under Rieffel-Morita
equivalence,
the C*-algebra of compact endomorphisms of the given Hilbert
C*-module has stable rank one---and so each element of
the algebra of compact endomorphisms with unit adjoined
can be approximated in norm by an invertible
element of this C*-algebra---of norm one if the given element
is of norm one.
The final step is to note that, if
an element of the algebra of compact endomorphisms
of a Hilbert C*-module with unit adjoined (or even just
an 
adjointable endomorphism) approximately
preserves inner products on a given finite subset, and
in addition has norm one, then its absolute value
is close to the identity on these elements---as, for an
adjointable endomorphism $x$ of norm at most one and a
Hilbert C*-module element $\xi$,
$$\eqalign{
\langle (1-|x|) \xi, (1-|x|)\xi\rangle &=
\langle \xi, (1-|x|)^2\xi\rangle \cr
&\le \langle\xi, (1-|x|)(1+|x|)\xi \rangle\cr
&= \langle\xi, (1-x^*x)\xi\rangle\cr
&=\langle \xi, \xi\rangle - \langle x\xi, x\xi\rangle\cr}$$
---and so if in addition this
endomorphism is invertible then it is close, on the given
finite subset, to its unitary part---a unitary element
of the C*-algebra of compact endomorphisms with unit adjoined, and
so an inner automorphism as defined above.

Finally, let $X$ and $Y$ be countably generated
Hilbert C*-modules
such that $[X]\le [Y]$, and let us show that $X$ is
isomorphic to a subobject of $Y$. One has a 
diagram as above but without the upwards arrows.
As above, the downwards arrows may be modified
one by one by composing with inner automorphisms
in such a way as to ensure that
each square is arbitrarily close to commuting on
an arbitrary finite set. As in 2.2 of [6],
these finite sets may be chosen in such a way that
the diagram is approximately commutative in the sense
of 2.1 of [6], and hence by 2.1 of [6] there is
a (unique) contraction from $X$ to  $Y$
such that the diagram remains commutative---and this map
is easily seen to be an isomorphism of $X$ onto a 
subobject of $Y$.
\smallskip

{\bf 4.}\
The following consequence of 
Theorems 1, 2, and 3 (and in particular of Lemma 2.2 of [17])
is of interest from the point of view of the classification
of C*-algebras. (By [9], it implies that
arbitrary simple C*-algebras stably isomorphic
to a separable simple AI algebra are classified by
their Elliott invariant.)
\smallskip

{\bf{Corollary.}}\ {\it {Let $A$ be the inductive limit of
a sequence of C*-algebras $A_1\to A_2\to \cdots$, and let
$B$ be a hereditary sub-C*-algebra of $A$. It follows
that for every finite subset of $B$ there is a
sub-C*-algebra of $B$ approximately containing this finite subset and
isomorphic to the image in $A$ of a hereditary
sub-C*-algebra of $A_n$ for some $n$. 
Hence in particular,
if $A$ is separable and if either $A$ has stable rank one, or 
quotients of hereditary sub-C*-algebras of
$A_n$ are weakly semiprojective for each $n=1, 2, \cdots,$ then
$B$ is isomorphic to the inductive limit of a sequence
of such C*-algebras.}}
\smallskip

{\it{Proof.}}\ 
The proof of the first assertion is very similar to part of the
proof of Theorem 2, above, but is self-contained---in the sense
that it appeals directly (in the same
way as above) to Lemma 2.2 of [17] and its proof.

Given a finite subset $F$ of $B$, choose a positive
element $h$ of $B$ of norm one which is close to the
identity on $F$, acting in either side, and choose 
a positive element $h_m$ of the image of $A_m$
in $A$ for some $n$ that is close to $h$, say to strictly
within distance $\epsilon$. Then by Lemma 2.2 of [17],
$$(h_n-\epsilon)_+ = dhd^*$$
for some $d\in A$
of norm at most one. (Here we need the contraction property---when this theorem
is used earlier we could just as well have used Lemma 2.5(ii) of [16] in which
the contraction property is not assured.)
Furthermore, inspection of the
construction of $d$ in [17] shows that, if $h_m$ is
sufficiently close to $h$, say to within distance
$\epsilon'$ (where $\epsilon'\leq  \epsilon$), then the
element $h^{1\over 2} d^* d h^{1\over 2}$
of $A$  is also close to $h$. Then with $a=dh^{1\over 2}$,
the hereditary sub-C*-algebras of $A$ generated by
$a^* a=h^{1\over 2} d^* dh^{1\over 2}$ and
$aa^*=dhd^* =(h_m-\epsilon)_+$ are
isomorphic. The first, that generated by $a^* a$, is contained 
in $B$ and
contains the finite set $a^* aF a^*a$. Since
$a^* a$ is close to $h$ and is of norm at most
one (as $h$ and $d$ are),
so that $a^* a F a^* a$ is close to $h F h$, and since
$h F h$ is close to $F$ (as $h$ is close to
the identity on $F$ on either side and has norm one),
the finite set $a^* aF a^*a$ is close to $F$. The
hereditary sub-C*-algebra of $A$ generated by $a^* a$
thus is contained in $B$ and approximately contains $F$.
The second hereditary sub-C*-algebra, that generated by $aa^*$,
is isomorphic to this, and, since $a^*a=(h_n-\epsilon)_+$,
for $m$ sufficiently large the
pre-image of the finite subset of the first algebra approximating
$F$ is approximated by the closure of $(h_m-\epsilon)_+
A_n (h_m-\epsilon)_+$,
where $A_n$ denotes also the image of $A_n$ in $A$.
This closure is then the image in $A$ of a hereditary
sub-C*-algebra of $A_n$, and, as a subalgebra
(although not a hereditary subalgebra)
of $A$ is isomorphic to a sub-C*-algebra of
$B$ approximately containing $F$, as desired.

Consider now the second assertion of the theorem.

If the sub-C*-algebras of $B$ constructed in the
way described above are weakly semiprojective,
as they are if all quotients of hereditary sub-C*-algebras
of $A_1, A_2, \cdots$ are assumed to be weakly semiprojective,
and if $B$ is separable, then an iterative 
construction using weak semiprojectivity together with an
interwining argument yields an increasing sequence of such
subalgebras with union dense in $B$, as desired.

An outline of this construction follows. (The proof is
slightly indirect.)
Choose first a sequence $D_1, D_2, \cdots$ of such subalgebras
such that for each $i$, a suitably large
subset of $D_i$ (to be specified) is approximately
contained in $D_{i+1}$. Then, by weak semiprojectivity
(maps from $D_i$ into the asymptotic sequence algebra
$\Pi^\infty_{n=i} D_n/\bigoplus^\infty_{n=i} D_n$ lift), we 
obtain, after passing to a subsequence of $D_1, D_2, \cdots,$
a sequence of C*-algebra maps $D_1\to D_2\to \cdots$ such
that the diagram
\smallskip
$$D_1 \to D_2 \to \cdots$$
\vskip 0.0051truein
$$\hskip -.5truein \downarrow\qquad \downarrow$$
\vskip 0.0051truein
$$B\,  \to\, B\ \to \cdots$$
approximately commutes, in the sense of 2.2 and 2.1 of [6],
where the maps $B\to B$ are the identity and the
vertical maps are the inclusions. By 2.1 of [6], the
diagram extends to an approximately commutative
diagram including a (unique) map 
$\lim\limits_\to D_i \to B=\lim\limits_\to B$, and
if the subalgebras $D_1, D_2,\cdots$ of $B$ are chosen to
approximate a dense sequence of elements of $B$ in a
suitable way (finitely many at each stage, eventually more and 
more, better and better), then the map $\lim\limits_\to D_i\to B$ is
surjective, as desired. (Cf.~[23].)

Now suppose that $A$ has stable rank one. In this
case, we can still prove the second assertion of the theorem,
but only by applying the full force of
Theorems 1, 2, and 3. (Not just as a consequence
of the first assertion---at least not
as simple a one as in the case considered above---the
intertwining argument for which was already somewhat 
indirect!)
We, in fact also need to use the proof of Theorem 2, not
just the statement---and not just the description of
inductive limits in the category $\cC u$ and 
how the functor from the category of C*-algebras to $\cC u$ 
acts on maps, but in fact the detailed construction of
how an element $x$ of $\cC u (A)$ is expressed as the
supremum of an increasing sequence
$$x_1 \leq x_2\leq \cdots$$
with $x_i$  the image of an element of $\cC u (A_i)$. 

Since, now, $B$ is separable, $B$ is singly generated as a hereditary 
sub-C*-algebra of $A$.
(In fact, for the present
case of the second assertion---that $A$ has stable
rank one---it is sufficient only to assume that $B$ is
singly generated---not necessarily separable as a
C*-algebra.)

The closed right ideal generated by $B$ is then
a countably (in fact singly) generated right
Hilbert $A$-module; denote this by $X$. By Theorem 2,
together with the concrete description of the inductive
limit in the category $\cC u$ of the
sequence 
$\cC u(A_1)\rightarrow \cC u(A_2)\rightarrow\cdots$ 
corresponding to the sequence of C*-algebras $A_1\to A_2\to\cdots$, given in the proof
of Theorem 2, there exists an increasing sequence
$x_1\leq x_2\leq \cdots$ in $\cC u (A)$ with $[X]=\sup x_i$, such that, for each $i$, the
element $x_i$ arises from a Hilbert $A_i$-module, say $X_i$
(by means of the natural map from $\cC u (A_i)$ to $\cC u (A)$), and
such that, moreover, $[X_i]\leq [X_{i+1}]$
in $\cC u(A_{i+1})$ for each $i$. In fact, as 
inspection of the construction in the proof of Theorem 2 shows, the
Hilbert C*-modules $X_1,\, X_2,\cdots$ over, respectively,
$A_1, A_2,\cdots$ may be chosen such that the push-forward $(X_i)_{A_{i+1}}$ is
isomorphic to a subobject of $X_{i+1}$ (over $A_{i+1}$) for
each $i$.

Note that, in particular, one has a sequence of
Hilbert $A$-module mappings (preserving the $A$-valued inner product)
$$(X_1)_A \to (X_2)_A \to \cdots,$$
and as shown in the proof of Theorem 1 it follows that
$$\sup x_i=[\lim\limits_\to (X_i)_A].$$
Since $x=[X]$ it follows that
$$[X]=[\lim\limits_\to (X_i)_A].$$
Since $A$ has stable rank one, by Theorem 3 we
have
$$X\cong \lim\limits_\to (X_i)_A,$$
where the isomorphism is as Hilbert $A$-modules. In particular,
the C*-algebra of compact
endomorphisms of $X$, i.e., $B$, is isomorphic to the
inductive limit of the sequence
$$B_1 \to B_2 \to \cdots$$
where $B_i$ denotes the C*-algebra of compact endomorphisms
of $(X_i)_A$ and $B_i\to B_{i+1}$ the canonical extension map
for compact endonorphisms. 
Now recall that, for each $i$, the Hilbert $A$-module $(X_i)_A$ arises from the
Hilbert $A_i$-module $X_i$, and that $(X_i)_{A_{i+1}}$ is
isomorphic to a subobject of $X_{i+1}$ (over $A_{i+1}$). Note
that, for each $i$, the C*-algebra $B_i$ of compact endomorphisms
of $(X_i)_A$ is the inductive limit of the natural sequence
of C*-algebras of compact endomorphisms of $(X_i)_{A_i}$,
$(X_i)_{A_{i+1}},\cdots$. Let us denote these C*-algebras by
$B^i_i, B^{i+1}_i, \cdots$. Thus, for each $i$,
$$B_i =\lim\limits_\to B^j_i$$
(where the limit is over $j$ with $j\ge i$). The preceding
statement implies that, for each $i$, the C*-algebra
$B^i_i$ of compact endomorphisms of $(X_i)_{A_i}$ is mapped
into the C*-algebra $B^{i+1}_{i+1}$ of compact
endomorphisms of $(X_{i+1})_{A_{i+1}}$, by the combination of
the canonical push-forward to the algebra $B^{i+1}_i$ of compact
endomorphisms of $(X_i)_{A_{i+1}}$ combined with the
canonical map from $B^{i+1}_i$ to $B^{i+1}_{i+1}$ corresponding
to the canonical extension of compact endomorphisms of 
$(X_i)_{A_{i+1}}$ to compact endomorphisms 
of $(X_{i+1})_{A_{i+1}}$, with
respect to the given isomorphism of $(X_i)_{A_{i+1}}$
with a subobject of $(X_{i+1})_{A_{i+1}}$.
By compatibility of the resulting infinite triangular
commutative diagram $(B^j_i)_{j\ge i}$ with the horizontal diagram 
$B_1\to B_2\to \cdots \to B$
it follows that
$$B=\lim\limits_\to B^i_i,$$
in fact as desired since by inspection of the construction
each $X_i$ has the same number of generators as $X$,
namely, one, and so $X_i$ is isomorphic to a closed right ideal
of $A_i$, and so $B^i_i$ to a hereditary 
sub-C*-algebra of $A_i$.
\smallskip

{\bf 5.}  The following consequence of Theorems 1, 2, and 3 is 
also of interest from the point of view of the classification
of C*-algebras.  Taken together with the preceding result, it might
be viewed as an indication of the potential usefulness of the Cuntz
invariant for proving isomorphism of C*-algebras.

\smallskip
{\bf Corollary.} {\it Let $A$ be a C*-algebra of stable rank
one.  An element of the Cuntz semigroup is compactly contained in
itself (in the order-theoretic sense)---let us refer to such an
element as compact---if, and only if, it corresponds to a Hilbert
C*-module which is algebraically finitely generated and projective.
(Hence in this case any Cuntz equivalence class also has these
properties.)  

Furthermore, $A$ is of real rank zero (every hereditary
sub-C*-algebra has an approximate unit consisting of projections---see
[3]) if, and only if, in the Cuntz semigroup of $A$, every element is
the supremum of an increasing sequence of compact elements.  
Alternatively,
an equivalent condition on the Cuntz semigroup is that an element
$x$ is compactly contained in an element $y$ (i.e., $x << y$) 
exactly when $x \le z \le y$ for some compact element $z$ (in other
words, $x \le z << z \le y$ for some $z$). 
(A different characterization of real rank zero, in the special case that
the Cuntz semigroup is almost unperforated, was given in [22].)

In particular, if $A$ is a separable, simple, AH algebra of stable 
rank one (for instance, by [11], if $A$ has an AH inductive
limit decomposition with diagonal maps between building blocks),
and if $A$ has the same Cuntz semigroup as a real rank zero AH
algebra with an AH inductive limit decomposition with building
blocks with spectra of bounded finite dimension, then by [19] and [20] (see also [18])
 $A$ is also such an
algebra. (Hence by [8], if also K$_1 (A) \cong$ K$_1 (B)$ then $A \cong 
B.)$}

\smallskip
{\it Proof}.  Let $X$ be a Hilbert C*-module over $A$ which is, just 
considered algebraically as a module, finitely generated
and projective.  As we shall show below, this is equivalent to the
property that $X$ is isomorphic as a (right) Hilbert $A$-module to 
the submodule of the finite direct sum $A^{\smallsim} \oplus \dots
\oplus A^{\smallsim}$ of finitely many
copies of $A$ determined by a projection $e$ in the
compact endomorphism C*-algebra of the Hilbert 
$A^{\smallsim}$-module $A^{\smallsim} \oplus\cdots \oplus
A^{\smallsim}$---i.e., to the submodule
$e(A^{\smallsim} \oplus\cdots\oplus A^{\smallsim})$,
where $A^{\smallsim}=A+\IC 1$ is the C*-algebra with unit adjoined,
considered as a module over $A$ in the natural way---with the entries
of $e$ assumed to belong to $A$, so that the $A^{\smallsim}$-valued
inner product
on the submodule determined by $e$ takes values in $A$---so
that it is indeed 
a Hilbert $A$-module.
In particular, one sees that the
identity endomorphism of $X$ is compact,
as the compact endomorphisms of $X$ are just the
restrictions to $X$ of the compact endomorphisms of
$A^{\smallsim}\oplus \cdots\oplus A^{\smallsim}$ 
taking $X$ into $X$---and these can be
identified with the unital algebra $eEe$ where $E$ denotes
the algebra of compact endomorphisms of
$A^{\smallsim}\oplus\cdots\oplus A^{^{\frak{\sim}}}$.

To check the assertion above, note just that,
purely algebraically, as an $A$-module, $X$ must be
isomorphic to $e(A^{\smallsim}\oplus\cdots\oplus A^{\smallsim})$
for some idempotent
element $e$ of M$_n(A^{\smallsim}$), where $n$ is the 
number of copies of $A^{\smallsim}$ in the direct sum (and can be
taken to be the number of generators of the module).
This is standard if $A$ is unital, but in the general case
one can just adjoin a unit to $A$ and note that modules over
$A$ are in bijective correspondence with unital modules over
$A^{\smallsim}$.
Next, recall that by Theorem 26 of [13], $e$ is similar in 
M$_n(A^{\smallsim})$ to a self-adjoint idempotent (i.e., projection),
and so we may suppose that $e$ is self-adjoint. Finally,
note that, presumably, any two Hilbert
$A$-modules which are isomorphic as $A$-modules are
in fact isomorphic as Hilbert $A$-modules, but, this
is quite elementary in the case that one of them is
$e(A^{\smallsim}\oplus\cdots\oplus
A^{\smallsim})$ as above. Namely, any module map from
$e(A^{\smallsim}\oplus\cdots\oplus A^{\smallsim})$ to a
Hilbert $A^{\smallsim}$-module (or even
just to a Banach $A$-module) $X$ must be continuous,
since this is trivially true in the case
$e=1$---and just $A^{\smallsim}$ in place of
$A^{\smallsim}\oplus\cdots\oplus A^{\smallsim}$.
In fact, it must be compact, since the identity of
$e(A^{\smallsim}\oplus\cdots\oplus A^{\smallsim})$ is
compact. (Since this is even finite rank, one sees that the composed map is also
of finite rank---and so we do not have to prove it is continuous.)
Then just pass to the unitary part of the polar decomposition of this continous
invertible map---to obtain an isomorphism of Hilbert $A^{\smallsim}$-modules,
from the existence of which one concludes that
$e(A^\sim\oplus\cdots\oplus A^{\smallsim})$
is in fact also a Hilbert A-module (i.e., the $A^{\smallsim}$-valued
inner product
takes values in $A$). (That we are in the
setting of elementary C*-algebra theory may be seen by noting that the map belongs to
the compact endomorphism C*-algebra of the direct sum of the two modules.)

Conversely, let $X$ be a countably generated Hilbert
$A$-module which, in 
the Cuntz semigroup, is order-theoretically compact.
Let us show that the C*-algebra of compact endomorphisms
of $X$ is unital. If not, since it has a strictly positive element
(see above), it has an  increasing
approximate unit $0\le u_1\le u_2\le \cdots$ such that
$u_{i+1} u_i=u_i$
and $u_{i+1}\not=u_i$ for every $i$. As shown earlier,
$$[X]=\sup[(u_i X)^-],$$
and here by compactness, $[X]\le [(u_i X)^-]$ for some $i$.
By Theorem 3, as $A$ has stable rank one $X$ is
isomorphic to a subobject of $(u_i X)^-$. It follows
that if $X\not=(u_i X)^-$ then the C*-algebra $B$
of compact endomorphisms of $X$ has a scaling element
in the sense of [2].

(If $b$ is a self-adjoint element of $B$ equal to the identity on $X_i=(u_iX)^{-}$,
and $v$ is an isomorphism from $X$ to a subobject of $X_i$, then $vb\in B$ (indeed,
to see this it is enough to consider the case $b$ has finite rank, say $b$ is the map
$\xi\mapsto \eta \langle \zeta,\xi \rangle$ for some $\eta,\zeta\in X$,
and then $vb$ is the map
$\xi\mapsto v\eta \langle \zeta,\xi \rangle$). With $x=vb$, we have
$x^*x=bv^*vb=b^2$ and $xx^*=vb^2v^*$, and so $(x^*x)(xx^*)=b^2vb^2v^*=vb^2v^*=xx^*$. 
Furthermore,
if $x^*x=xx^*$ then $b^4=b^2$, i.e., $b^2$ is a projection, and, furthermore,
$vb^2v^*=b^2$, which implies that $b^2$ is exactly the unit of the compact endomorphism algebra
of $X_i$, which is impossible by hypothesis. Thus, $x^*x\neq xx^*$, and so $x$ 
is a scaling 
element.)  

Hence, by the proof of Theorem 3.1 of [2], the C*-algebra 
$A^{\smallsim}$ obtained by adjoining a
unit to $A$ has an infinite projection.
This contravenes the hypothesis that $A$---or rather (by definition),
$A^{\smallsim}$---has
stable rank one (used now for the second time
in the proof of this implication).
(To show that a unital C*-algebra of stable rank one does not have an
infinite projection, it is enough to
show that the unit cannot be infinite, or, equivalently,
that every isometry is unitary. But if an isometry can be
approximated arbitrarily closely by an invertible element, then 
it can also be approximated arbitrarily closely by the unitary
part of this invertible element, and so its
range projection also can be approximated arbitrarily closely
by the unit, and so must be the unit; the isometry
must be unitary.) 

Now consider the second assertion of the theorem. Suppose first that $A$ has real rank zero,
and let us show that every element of $\cC u(A)$ is the supremum of an increasing 
sequence of compact elements.  Let $X$ be a countably generated Hilbert $A$-module. Since
by Theorem 2 of [14], $X$ is isomorphic to a closed submodule of 
$A\oplus A\oplus\cdots$ (cf.~just above),
the C*-algebra of compact endomorphisms of $X$ is isomorphic to a hereditary sub-C*-algebra 
of the 
stabilization of $A$ and is therefore also of real rank zero. In particular---since it has a
countable approximate unit (see proof of Theorem 1)---this C*-algebra has an increasing
sequential approximate unit $(e_n)$ consisting of projections. Then
$X=\lim\limits_\to e_nX$, and so by the proof of Theorem 1,
$[X]=\sup [e_nX]$ in $\cC u(A)$.  As shown above, $[e_nX]$ is compact in $\cC u(A)$ for each $n$.
(The proof of this first implication does not use that $A$ is simple or has stable rank one.)

Now suppose that every element of $\cC u(A)$ is the supremum of an increasing sequence of compact elements, and let us show that $A$
has real rank zero. By Theorem 2.6 of [3], 
it suffices to show that any singly generated hereditary sub-C*-algebra
of $A$, say $B$, has an increasing approximate unit consisting of projections, or, equivalently,
that the closed right ideal $B$ generates, say $X$, is the closure of an increasing sequence of closed
right ideals generated by projections,---or, equivalently again, that $X$ is the Hilbert 
$A$-module inductive
limit of a sequence of closed right ideals of $A$ generated by projections, with respect to
Hilbert $A$-module isomorphisms between the objects at various stages and subobjects of the
succeeding ones. (This is clearly equivalent to the same statement for the increasing
sequence of images of these objects in $X$, and any subobject of $X$ is a closed right
ideal, necessarily generated by a projection if this is true up to isomorphism, by the first
assertion of the theorem.)

By hypothesis, $[X]=\sup\, x_i$ in $\cC u(A)$ where $x_1\leq x_2\cdots $ is an increasing
sequence of compact elements of $\cC u(A)$. With $X_1, \, X_2\,\,\cdots$ Hilbert C*-modules such
that $x_i=[X_i]$ for each $i$, by Theorem 3 (as $A$ has stable rank one) $X_i$ is isomorphic
to a subobject of $X_{i+1}$ for each $i$. As shown in the proof of Theorem 1, in this case
$$\sup\,[X_i]=[\lim_{\rightarrow} X_i].$$
By Theorem 3 again (actually, the first use of Theorem 3 above was not necessary on account of the 
very special nature of the Hilbert $A$-modules $X_i$---algebraically finitely generated
and projective by the first assertion of the theorem),
$$X\cong \lim\limits_\to X_i$$
where $\cong$ denotes isomorphism of Hilbert $A$-modules. In particular, we may suppose
that each $X_i$ is a subobject of $X$ and that the maps $X_i\to X_{i+1}$ are inclusions.
Since $X$ is a closed right ideal of $A$ it follows that $X_i$ is also for each $i$,
and it remains to note that also $X_i=e_iA$ for a projection $e_i\in A$. (Otherwise $X_i$
could not be finitely generated algebraically, let alone be projective!)
\smallskip



{\bf{6. Appendix.}}\ Let us explain in more detail the relationship
between the ordered semigroup $\cCu A$  
defined in Section 1 (and studied in Theorems 1, 2, 
and 3 and Corollaries 4 and
5 above) and the ordered semigroup 
introduced by Cuntz in [4]---denoted by W$A$ by R\o rdam
in [25] and now often called the Cuntz semigroup. 
(In [25] the structure of W$A$ was referred to as that
of positive ordered semigroup,
to reflect the fact that every element is positive. In the present 
article we suppress this additional qualifier, since we are only
considering such ordered abelian semigroups.)
(The main purpose of our proposed new notation is to emphasize 
the additional structure we have identified---the operation 
of taking countable increasing suprema
and the relation of compact containment in the order-theoretic 
sense.) 

Briefly, if $A$ is stable, i.e., if $A$ is isomorphic to 
$A\otimes  \K$ 
where $\K$  denotes the C*-algebra of compact operators on a 
separable infinite-dimensional Hilbert space, then the two 
ordered semigroups are exactly the same (the functors
are equivalent by a natural transformation).
 
 In general, $\cCu A$ is the same as $\cCu (A\otimes \K$) 
(and so the same as W($A\otimes \K$))

Let us show first that, with $\cCu A$ defined as in Section 1,
$\cCu A$ is isomorphic to $\cCu (A\otimes \K$)---and that the
isomorphism may be chosen to be natural (i.e.,
to be a natural transformation between these two functors).

By the functoriality of $\cCu$ (Theorem 2 above), corresponding
to the inclusion of $A$ as $A\otimes  e$ in $A\otimes \K$,
where $e$ is a (fixed) minimal non-zero projection in $\K$, there is a
morphism $\cCu A \to \cCu (A\otimes \K$) in the category $\cCu$, which 
furthermore (again by functoriality) constitutes a natural
transformation between these two functors. It remains to show 
that this map is an isomorphism for any C*-algebra $A$. 

On specializing to the present case, the morphism
$\cCu A\to \cCu (A\otimes \K$) 
is seen to consist of, given a countably generated
(right) Hilbert $A$-module $X$, taking the 
Hilbert $A$-module direct sum of a countable infinity of copies
of $X$, and then having both $A$ and $\K$---and therefore also
$A\otimes \K$---act on this in the natural way (on the right). 
Let us show that the map in the opposite direction, beginning 
with a countably generated 
Hilbert $A\otimes \K$-module $Y$, and cutting it down by the
subalgebra $A\otimes e$ where $e$ is a fixed minimal non-zero 
projection in $\K$, resulting in a Hilbert C*-module over this 
C*-algebra which is naturally isomorphic to $A$, 
and which is countably generated if $Y$ is, preserves our
notion of equivalence of countably generated Hilbert C*-modules, 
and at the level of Cuntz semigroups (in the sense of the 
present paper) is the inverse of the
map just described. To see that equivalence is preserved by
the backwards map, it is enough to note that this map 
preserves isomorphism (it in fact
preserves arbitrary homomorphisms in a natural way), and both
this map and the map in the forwards direction preserves the relation
of inclusion and the relation of compact containment (for a
subobject---see Theorem 1). It is straightforward that these 
maps are inverse to each other at the level of Hilbert 
C*-modules, and the desired isomorphism of the Cuntz 
semigroups follows.

Let us now show that, if $A$ is stable, then the map which to a
positive element of $A$ associates the closed right ideal it 
generates, considered as a Hilbert $A$-module, determines an
isomorphism between W$A$ and $\cCu A$.

Recall that two positive elements of $A$, let us say $a$ and $b$,
are comparable in the sense of Cuntz, with $b$ majorizing $a$,
if there exists a sequence
($c_n$) in $A$ such that $c_nbc_n^*$ converges to $a$. Let us
suppose that this holds, and let us show that the class in 
$\cCu A$ of the closure of $aA$ is
majorized by the class of the closure of $bA$.

By Lemma 2.2 of [17], for each $n$ there exists $d_n$ in $A$ such that
$d_nc_nbc_n^*d_n^*$ is a continuous function $a_n$ of $a$ and the sequence
($a_n$) is increasing with limit $a$. It follows on the one hand that
the closure of $a_nA$ is isomorphic to a subobject of the closure of
$bA$ for each $n$, and in particular the class of this Hilbert
$A$-module in $\cCu A$ is majorized by the class of the closure
of  $bA$, and on the other hand, as shown in the proof of Theorem 1, 
that the class of the closure of $aA$ in $\cCu A$ is the supremum
of the increasing sequence of classes of the closures of 
$a_1A,\, a_1A,\, \cdots$.
It follows immediately that the class of the closure of $aA$ 
in $\cCu A$ is majorized by the class of the closure of $bA$.

Let us show, conversely, that if $a$ and $b$ are as above, and the
closed right ideal of $A$ generated by $a$ is majorized in the ordered
semigroup $\cCu A$ defined above (actually proved to be an
ordered semigroup only in the proof of Theorem 1) by the closed
right ideal generated by $b$, then $a$ is majorized by $b$ in the sense
of Cuntz. Choose a continuous positive real-valued function $f$ on
the spectrum of $a$ equal to zero in a neighbourhood of zero (if
zero belongs to the spectrum of $a$), such that $f(a)$ is close to $a$.
Then the closed right ideal generated by $f(a)$ is compactly contained
in that generated by $a$ (as there exists a continuous function $g$
on the spectrum of $a$, equal to zero at zero, such that 
$g(a)f(a) = f(a))$,
and is therefore by hypothesis isomorphic, as a Hilbert $A$-module,
to a subobject (compactly contained, but we shall not need this)
of the closed right ideal generated by $b$. Such a subobject
must in fact be a smaller closed right ideal, countably generated
and therefore singly generated. 
It follows that there is an element $x$ of $A$ such that
$x^*x = f(a)$
and $xx^*$ generates the closed right ideal in question.
Since then $x$ is also in this ideal, which is contained in the
closed right ideal generated by $b$, there is  
a sequence ($a_n$) in $A$ such that $ba_n$ converges to $x$,
and then $f(a)$ is
the limit of $x^*ba_n$. By  
polarization (namely, the equation $cbd^*$ equal to the
average of the elements $(c+zd)b(c+zd)^*$ with $z$ a power of $i$,
which holds as $b$ is self-adjoint), $x^*ba_n$ is majorized by $b$
in the sense of Cuntz
for each $n$, and hence, since the set of such elements is closed
(as follows immediately from the definition),
also $f(a)$ is majorized by $b$ in this sense---indeed, since $f(a)$ is
arbitrarily close to $a$, also $a$ is.

It follows that the map from positive elements of $A$ to closed
right ideals gives rise to an isomorphism of ordered semigroups
between W$A$ and $\cCu A$, or, rather, 
in the first instance,
between the subset of W$A$ arising from
positive elements of $A$, as opposed to matrix algebras over $A$,
and the subset of $\cCu A$
arising from singly generated closed right ideals of $A$,
as opposed to countably
generated Hilbert C*-modules over $A$---but these subsets 
exhaust these two semigroups in the case that $A$ is stable,
as we shall now show. 

First, let a be a positive element of $A\otimes \rM_n$,
and let us show that it
is equivalent in W$A$ to a positive element of
$A\otimes e$, where $e$ is a non-zero
minimal projection in M$_n$, using of course that $A$ is
stable. Writing $A$ as
$B\otimes \K$ for some $B$, and noting that there is an
isometry $v$ in the multiplier
C*-algebra of $\K\otimes \rM_n$ such that
$v(\K\otimes \rM_n)v^* = \K\otimes e$, we have that
($1\otimes v)a(1\otimes v)^*$ is a positive element of
$\K\otimes e$, which is easily
seen to be equivalent in W$A$ to $a$.
(Recall that $v$ is the limit of a sequence of
elements of $\K\otimes \rM_n$) in the strict topology on the
multiplier algebra.)

Second, and finally, let $X$ be a countably generated
Hilbert C*-module over $A$, and let us show that $X$ is
isomorphic to a singly generated closed
right ideal of $A$ (using again that $A$ is stable).
It is enough to show that, as a Hilbert $A$-module, $X$ is
isomorphic to just some closed
right ideal of $A$, since this is then countably generated
as $X$ is, and
a countably generated closed right ideal of a C*-algebra is singly
generated. (A countable set of generators may be assumed to be
positive and summable in norm, and the sum is then a single 
generator---as the
self-adjoint part of a closed two-sided ideal is a hereditary
sub-C*-algebra.)

By Theorem 2 of [14],
$X$ is isomorphic to a direct summand of the Hilbert
C*-module direct sum of
a countable infinity of copies of the Hilbert $A$-module $A$,
and in particular
to a closed submodule of this direct sum. In fact, since
$A$ is stable, this
infinite direct sum Hilbert $A$-module is isomorphic to $A$!

(To see this, write $A$ as $B\otimes \K$ for some C*-algebra $B$,
and note that
the Hilbert $A$-module $A$ is equal to the (internal) 
Hilbert C*-module direct
sum of the closed right ideals ($1\otimes e_n)(B\otimes \K$)
where ($e_n$) is
a sequence of minimal closed two-sided projections in $\K$ 
with sum equal
to one in the multiplier C*-algebra of $\K$. In other words, $A$ 
is isomorphic
as a Hilbert $A$-module to an infinite direct sum of copies of some
Hilbert $A$-module, and partitioning the index set into a countable
infinity of subsets of the same cardinality as the whole set 
one sees that $A$ is isomorphic to a countably infinite direct 
sum of copies of itself.)

It follows that $X$ is isomorphic to a closed right ideal of $A$,
and a countably generated one since $X$ is countably generated. 
Recall, finally, that a countably generated closed right ideal of 
a C*-algebra is singly generated.  

   One advantage of the original description of the Cuntz
semigroup of a C*-algebra is that it is immediate that 
approximately inner automorphisms
of the C*-algebra act trivially on it. 
(In the Hilbert C*-module
setting, which is remarkably useful for a number of purposes, as may
be clear by now---for instance, in the case of a C*-algebra of stable
rank one Cuntz equivalence just amounts to isomorphism of
Hilbert modules, in perfect analogy with Murray-von Neumann 
equivalence---it is
perhaps not quite obvious that even inner automorphisms act
trivially. In fact, a straightforward algebraic calculation 
establishes this.)

Another advantage of the original description of the Cuntz semigroup
is that, for a non-stable C*-algebra, it contains additional
information: just as in the case of the Murray-von Neumann semigroup,
one may keep track of when additional classes appear when one passes
to matrix algebras---or, for that matter, when one stabilizes 
(although in the case of the Murray-von Neumann semigroup this 
last step introduces no new classes). On the other hand, 
this information is also readily
expressible in the Hilbert module language. Just as projections
in a matrix algebra of a certain order over an arbitrary algebra 
correspond to projective modules having a generating subset 
with that number of elements, so
also do Cuntz classes arising in the original sense from a
matrix algebra of a certain order over a C*-algebra correspond 
to Hilbert C*-modules
over the C*-algebra with that (finite) number of
generators---it is only when
one looks at Cuntz classes arising from the stabilization
that one obtains a
Hilbert C*-module requiring an infinite number of generators.
(On the other hand,
for most purposes it is already of interest to consider the case
of a stable C*-algebra, in which case every countably generated
Hilbert C*-module is singly generated, as a Hilbert C*-module.)
(Note that a Hilbert C*-module which is countably generated
purely algebraically, i.e., as a module, must in fact be 
finitely generated and 
projective, as a module.)

\bigskip
\bigskip

{\bf References}
\bigskip
\itemitem{1.\ }
B. Blackadar, {\it K-Theory for Operator Algebras (Second edition),}
Mathematical Sciences Research Institute Publications, {\bf 5}, Cambridge
University Press, Cambridge, 1998.
\smallskip

\itemitem{2.\ }
B. Blackadar and J. Cuntz, {\it The structure of stable algebraically
simple C*-algebras}, Amer. J. Math. {\bf 104} (1982), 813--822.
\smallskip

\itemitem{3.\ }
L. G. Brown and G. K. Pedersen, {\it C*-algebras of real rank zero}, J. Funct. Anal.
{\bf 99} (1982), 131--149.
\smallskip

\itemitem{4.\ } J. Cuntz, {\it Dimension functions on simple C*-algebras},
Math. Ann. {\bf 233} (1978), 145--153.
\smallskip

\itemitem{5.\ }
G. A. Elliott, {\it The ideal structure of the multiplier algebra of an AF algebra},
C. R. Math. Acad. Sci. Soc. R. Can. {\bf 9} (1987), 225--230.
\smallskip

\itemitem{6.\ }
G. A. Elliott, {\it On the classification of
C*-algebras of real rank zero},
J. Reine Angew. Math  {\bf 443} (1993), 179--219.
\smallskip

\itemitem{7.\ }
G. A. Elliott, {\it Towards a theory of classification}, preprint.
\smallskip

\itemitem{8.\ }
G. A. Elliott, G. Gong, and L. Li, {\it On the classification of
simple inductive limit C*-algebras  II: The isomorphism theorem,}
preprint, 1998 (to appear, Invent. Math.).
\smallskip

\itemitem{9.\ }
G. A. Elliott and C. Ivanescu, {\it The classification of separable simple
C*-algebras which are inductive limits of continuous-trace C*-algebras with spectrum
homeomorphic to the closed interval $[0,1]$}, preprint.
\smallskip

\itemitem{10.\ }
G. A. Elliott and K. Kawamura, {\it A Hilbert bundle characterization of Hilbert C*-modules},
preprint.
\smallskip

\itemitem{11.\ }
T. M. Ho, {\it On inductive limits of homogeneous C*-algebras with diagonal
morphisms between the building blocks}, Ph.D. thesis, University of Toronto, 2006.
\smallskip

\itemitem{12.\ }
K. K. Jensen and K. Thomsen, {\it Elements of KK-Theory}, Mathematics:
Theory and Applications, Birkh\"auser, Basel, 1991.
\smallskip

\itemitem{13.\ }
I. Kaplansky, {\it Rings of Operators}, W. A. Benjamin, Inc., New York, 1968.
\smallskip

\itemitem{14.\ } G. G. Kasparov, {\it Hilbert C*-modules: theorems of Stinespring
and Voiculescu}, J. Operator Theory {\bf 4} (1980), 133--150.
\smallskip

\itemitem{15.\ }
E. C. Lance, {\it Hilbert C*-modules. A tool kit for operator algebraists}, London
Mathematical Society Lecture Note Series, {\bf 210}, Cambridge University Press,
Cambridge, 1995.
\smallskip

\itemitem{16.\ }
E. Kirchberg and M. R\o rdam, {\it Non-simple purely infinite C*-algebras},
Amer. J. Math. {\bf 122} (2000), 637--666.
\smallskip

\itemitem{17.\ }
E. Kirchberg and M. R\o rdam, {\it Infinite non-simple C*-algebras absorbing
the Cuntz algebra $\Cal {O}_{\infty}$}, Adv. Math. {\bf 167} (2002), 195--264.
\smallskip

\itemitem{18.\ }
H. Lin, {\it An introduction to the classification of
amenable C*-algebras}, World Scientific Publishing Co., Inc.,
River Edge, NJ, 2001.
\smallskip

\itemitem{19.\ }
H. Lin, {\it Simple AH algebras of real rank zero}, Proc. Amer. Math. Soc. {\bf 131}
(2003), 3813--3819.
\smallskip

\itemitem{20.\ }
H. Lin, {\it Classification of simple C*-algebras of tracial topological rank zero},
Duke Math. J. {\bf 125} (2004), 91--119.
\smallskip

\itemitem{21.\ }
V. M. Manuilov and E. V. Troitsky, {\it Hilbert C*-modules. Translated from the 2001
Russian original by the authors.} Translations of Mathematical Monographs,
{\bf 226}, American Mathematical Society, Providence, RI, 2005.
\smallskip

\itemitem{22.\ }
M. R\o rdam, {\it On the structure of simple C*-algebras tensored with a UHF
algebra, II}, J. Funct. Anal. {\bf 107} (1992), 255--269.
\smallskip

\itemitem{23.\ }
M. R\o rdam, {\it A short proof of Elliott's theorem: $\Cal {O}_2\otimes \Cal{O}_2
\cong\Cal{O}_2$,} C. R. Math. Acad. Sci. Soc. R. Can. {\bf 16} (1994), 31--36.
\smallskip

\itemitem{24.\ }
M. R\o rdam, {\it Classification of nuclear, simple
C*-algebras,} pages 1--145 of {\it Classification of nuclear
C*-algebras. Entropy in operator algebras}.
Encyclopedia of the Mathematical Sciences, {\bf 126},
Springer, Berlin, 2002.
\smallskip

\itemitem{25.\ }
M. R\o rdam, {\it The stable and the real rank of $\Cal {Z}$-absorbing C*-algebras}, International J. 
Math. {\bf 15} (2004), 1065--1084.
\smallskip 

\itemitem{26.\ }
A. S. Toms, {\it On the classification problem for nuclear C*-algebras}, Ann. of Math.,
to appear.
\medskip

\vglue .3truein

Department of Mathematics,\par
University of Toronto,\par
Toronto, Ontario,\par
Canada\ \ M5S 2E4\par
\smallskip

kris\@math.toronto.edu\par
elliott\@math.toronto.edu\par
\medskip

Department of Mathematics,\par
Northern British Columbia,\par
Prince George, British Columbia,\par
Canada\ \ V2N 4Z9\par
\smallskip

cristian\@math.toronto.edu

\end